\documentclass{article}
\input epsf.tex
\usepackage{latexsym}
\usepackage{amsfonts,amssymb,amsmath}
\usepackage{epsfig}
\usepackage{psfig}   

\newtheorem{theorem}{Theorem}[section]
\newtheorem{lemma}[theorem]{Lemma}

\newtheorem{question}[theorem]{Question}

\newtheorem{conjecture}[theorem]{Conjecture}

\newtheorem{example}[theorem]{Example}
\newtheorem{convention}[theorem]{Convention}

\def\lk{\mbox{\rm{lk}}}
\def\Isom{\mbox{\rm{Isom}}}

\title{Small curvature surfaces in hyperbolic 3-manifolds}
\author{Christopher J. Leininger \thanks{This research was partially conducted by the author for the Clay Mathematics Institute and while supported by an N.S.F. postdoctoral fellowship.}}

\begin{document}

\maketitle

\begin{figure}[htb]
\begin{center}
\ \psfig{file=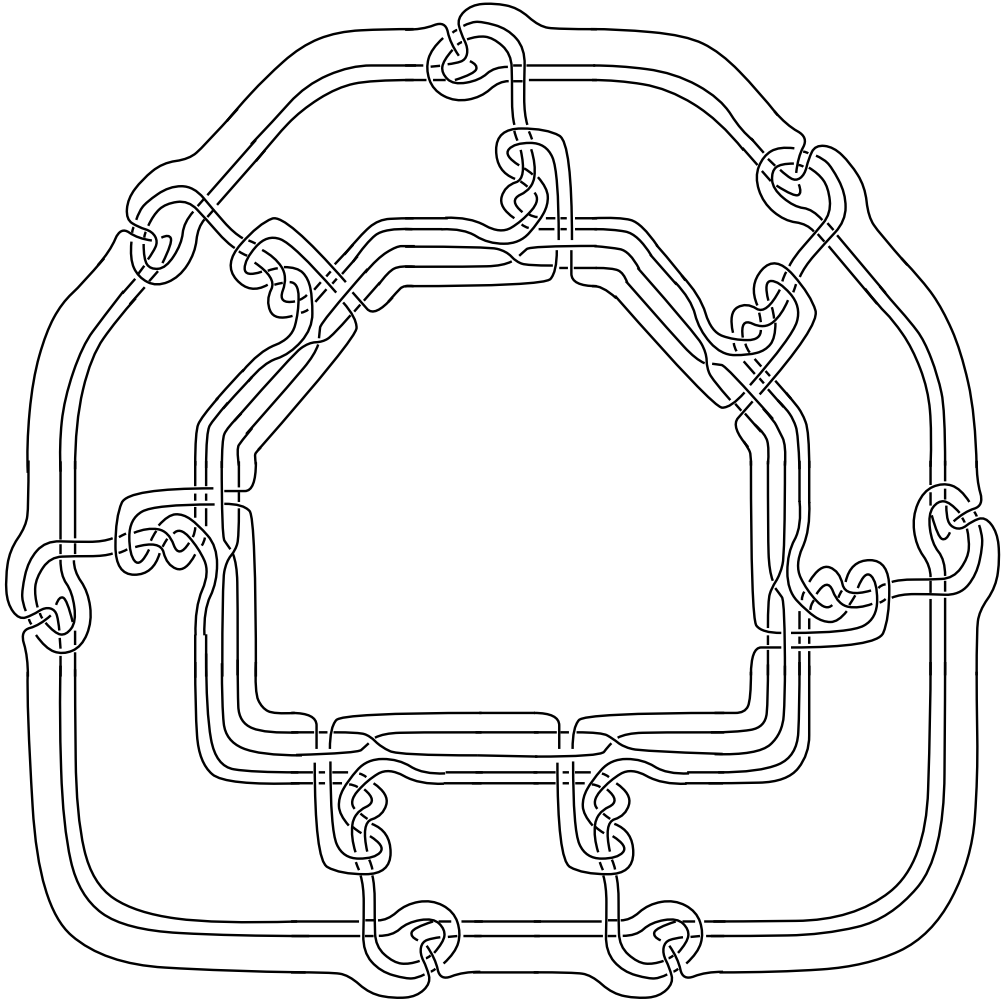,height=5.24truein}
\caption{A two component link containing a closed embedded totally geodesic surface of genus 6.}
\label{holycrap}
\end{center}
\end{figure}

\section{Introduction}

In the study of 3-manifolds, embedded essential surfaces have played a prominent role (see e.g. \cite{Waldhausen}, \cite{Tbull}, and \cite{CGLS}).
In \cite{MenascoReid}, Menasco and Reid study the geometry of closed embedded surfaces in hyperbolic link complements in $S^{3}$.
The main question under investigation is of existence of closed embedded totally geodesic surfaces in such manifolds.
In particular, they describe an eight component link in $S^{3}$ (whose construction they attribute to W. Neumann) with hyperbolic complement containing a closed embedded totally geodesic surface.
They also prove that the complement of a hyperbolic knot which is either alternating, 3-braid, or has 2-generator knot group cannot contain a closed embedded totally geodesic surface.
They make the following conjecture (see also \cite{Kirby}, Problem 1.76).

\begin{conjecture} \label{noknots}
There are no hyperbolic knots in $S^{3}$ which contain closed embedded totally geodesic surfaces in their complement.
\end{conjecture}

Since \cite{MenascoReid}, Conjecture \ref{noknots} has also been verified for toroidally alternating knots \cite{Adams2} (which includes almost alternating knots \cite{Adams_et_al} and Montesinos knots \cite{Oertel}), 3-bridge knots and double torus knots \cite{IchOz}, and $4$-braid knots \cite{Hiroshi}.

A totally geodesic surface $F$ in a hyperbolic 3-manifold $M$ with toroidal boundary has the property that $M / F$ is acylindrical (here $M / F$ denotes {\em $M$ cut along $F$}, which we take to be the path metric completion of $M \setminus F$).
Thus, failure of $M / F$ to be acylindrical is an obstruction to $F$ being totally geodesic.
In \cite{AdamsReid} it is shown that there exists hyperbolic knots $K \subset S^{3}$ with $S^{3} \setminus K$ containing surfaces $F$ such that $(S^{3} \setminus K ) / F$ is acylindrical.
This is essentially the only evidence for a counter-example to Conjecture \ref{noknots}.

In this paper we attempt to provide more evidence for such a counter-example, by showing that it is possible to get ``as close as possible'' to hyperbolic knots in $S^{3}$ with totally geodesic surfaces in their complements.
There are several meanings of ``close'' which we consider, and in all cases, we construct examples as close as possible in that sense.\\

To begin with, we could simply take ``close'' to mean that the number of components is as small as possible.
In Section \ref{linksection} we find links with only two components which contain embedded totally geodesic surfaces in their complements.\\

\noindent
{\bf Theorem \ref{links} }
{\em For any even integer $g \geq 2$, there exists a two component hyperbolic link in $S^{3}$ which contains an embedded totally geodesic surface of genus $g$ in its complement.}\\

Recall that for a surface in a Riemannian 3-manifold, the principal curvatures at a point measure the deviation of that surface from being totally geodesic at that point (see Section \ref{backgroundsection}).
Thus a hyperbolic knot with a surface of small principal curvature is ``close'' to a totally geodesic surface.
In section \ref{knotsection} we prove\\

\noindent
{\bf Theorem \ref{knots} }
{\em For any $g \geq 3$ and any $\epsilon > 0$, there exists a hyperbolic knot $K \subset S^{3}$ containing a closed embedded surface of genus $g$ in its complement whose principal curvatures are bounded in absolute value by $\epsilon$.}\\

We could also take ``close'' to mean that the knot complement has pinched negative sectional curvature with small pinching ratio, and it contains an embedded totally geodesic surface.
A slight modification of the construction of Theorem \ref{knots}, along with the ideas of the proof of the Gromov-Thurston $2\pi$-Theorem gives us\\

\noindent
{\bf Theorem \ref{pinchknots} }
{\em For any $g \geq 3$ and any $\epsilon > 0$, there exists a knot $K \subset S^{3}$ with complement supporting a Riemannian metric with negative sectional curvatures pinched between $-1-\epsilon$ and $-1+\epsilon$, which contains a closed embedded totally geodesic surface of genus $g$.}\\

Two other notions for getting ``close'' are known to hold.
First, one can remove the requirement that the surface be embedded.
In this case, we see from \cite{Maclachlan} and \cite{MacReid} that the figure eight knot complement contains infinitely many closed, {\em immersed}, totally geodesic surfaces.
Second, we can require only that our surfaces have finite area.
Students in an R.E.U. under the direction of Colin Adams have recently constructed knots in $S^{3}$ which contain embedded totally geodesic {\em cusped} surfaces (see \cite{Adams3}, \cite{Adams4}, and also \S \ref{knotcuspsection}).
For the sake of completeness, we have included such examples (Example \ref{cuspedsurfaceexample}) in \S \ref{knotcuspsection}.\\

We also note that not only do small principal curvature surfaces behave geometrically like totally geodesic surfaces, but they also exhibit many of the same {\em topological} properties.
Evidence of this was provided by Thurston who noticed that if the principal curvatures are strictly bounded by $1$ in absolute value, then the surface is incompressible (see \cite{AitchRub} and also Section \ref{smallcurvesurface}).
In Section \ref{smallcurvesurface} we continue with this comparison and show that the manifold obtained by cutting open along a surface with sufficiently small principal curvatures is acylindrical, with an obvious exception.
In particular, we prove\\

\noindent
{\bf Theorem \ref{curvature_acylindrical} }
{\em  Given $g, r > 0$, there exists $\delta > 0$ with the following property.
If $F$ is an embedded, closed, orientable surface in an oriented finite volume hyperbolic 3-manifold $M$, with $genus(F) \leq g$, $injrad(F) \geq r$, and principal curvatures of $F$ bounded by $\delta$ in absolute value, then either $F$ bounds a twisted $I$-bundle in $M$, or $M/ F$ is acylindrical.}\\

The possibility that $F$ bounds a twisted $I$-bundle is necessary since as $\epsilon \rightarrow 0$, the $\epsilon$-neighborhood of a non-orientable totally geodesic surface has principal curvatures approaching $0$.

We have presented Theorems \ref{links}, \ref{knots}, and \ref{pinchknots} as evidence for a negative resolution to Conjecture \ref{noknots}.
On the other hand, if Conjecture \ref{noknots} is true these theorems indicate the difficulty in a geometric approach to proving it.
As Theorem \ref{curvature_acylindrical} shows, it can be difficult to distinguish, both geometrically and topologically, between totally geodesic surfaces and surfaces with very small principal curvatures.
Furthermore, the totally geodesic property for a surface is unstable, and many of the coarse geometric methods for studying hyperbolic 3-manifolds are often a bit too insensitive to this.

We remark that from a number theoretic point of view (in terms of the holonomy representation of the 3-manifold), totally geodesic surfaces are easily distinguished from every other type of surface.
Moreover, arithmetic information, along with the topological information of being a knot complement in $S^{3}$ has already proven to be very restrictive indeed (see \cite{Reid}).\\

The paper is organized as follows:
Section \ref{backgroundsection} contains a few definitions and theorems from 3-manifold topology and Riemannian and hyperbolic geometry necessary for our work.
In Section \ref{linksection}, we construct the required hyperbolic links with totally geodesic surfaces in their complements.
Section \ref{knotsection} contains the various constructions of the required knots.
In Section \ref{smallcurvesurface} we prove Theorem \ref{curvature_acylindrical}.
Section \ref{relatedquestions} includes two questions related to Conjecture \ref{noknots}.\\

\noindent
{\bf Acknowledgments.} Thanks to Alan Reid, Lewis Bowen, and Joe Masters for helpful conversations regarding this work and to The University of Texas at Austin for allowing me to work there during the summer 2002.
Special thanks to Alan for his encouragement in writing this up.

\section{Background} \label{backgroundsection}

\subsection{$3$-manifolds}

Here we collect some of the basic facts and definitions concerning $3$-manifolds and orbifolds, and Dehn filling, see \cite{Hempel}, \cite{Rolfsen}, \cite{NeumannReid}, and \cite{Tnotes} for more details.

An orientable properly embedded surface $(F,\partial F)$ in a compact, orientable $3$-manifold (or $3$-orbifold) $(M,\partial M)$ is {\em incompressible} if $F \not \cong D^{2}$, $F \not \cong S^{2}$, and the inclusion induces an injection on (orbifold) fundamental group.
$F$ is {\em essential} if $F$ is incompressible and is not properly homotopic into $\partial M$.
$M$ is acylindrical if it does not contain an essential annulus.

Let $\partial_{0} M \cong T^{2}$ be a torus boundary component of $M$ (which we assume is disjoint from the singular locus for orbifolds).
A {\em slope} on $\partial_{0} M$ is an isotopy class of unoriented essential simple closed curves on $\partial_{0} M$.
Slopes on $\partial_{0} M$ are in a $2$ to $1$ correspondence with primitive elements of $H_{1}(\partial_{0} M;{\mathbb Z}) \cong {\mathbb Z}^{2}$, with the ambiguity coming from the lack of orientations.
We will often ignore this ambiguity, making no distinction between slopes and primitive elements of $H_{1}(\partial_{0} M;{\mathbb Z})$.

Let $\mu$ and $\lambda$ denote generators for $H_{1}(\partial_{0} M; {\mathbb Z})$.
Slopes on $\partial_{0} M$ then correspond to co-prime integer pairs in ${\mathbb Z}^{2}$ by associating the pair $(p,q)$ to the slope $p \mu + q \lambda$ (note that $(p,q)$ and $(-p,-q)$ represent the same slope).

Given a slope $\alpha$ on $\partial_{0} M$, one can form a new orbifold (or manifold, when $M$ is a manifold) $M(\alpha)$ by {\em $\alpha$-Dehn filling} on $\partial_{0} M$, as follows.
Let $S^{1} \times D^{2}$ be a solid torus.
Choosing a homeomorphism $h : \partial (S^{1} \times D^{2}) \rightarrow \partial_{0} M$, so that $h(* \times \partial D^{2})$ represents $\alpha$, we can glue $S^{1} \times D^{2}$ to $M$ by identifying points $x$ and $h(x)$.
The resulting space is an orbifold (or manifold), and up to homeomorphism, depends only on $\alpha$.
If we have chosen a basis $\mu, \lambda$ for $H_{1}(\partial_{0} M; {\mathbb Z})$ and $\alpha$ is given by $(p,q)$, we denote $M(\alpha)$ by $M(p,q)$.
Note that there is a natural inclusion $i : M \rightarrow M(\alpha)$.

A variation of this construction that we will make use of is {\em orbifold Dehn filling}, which we now describe.
Given an integer $d > 1$ and a slope $\alpha$ on $\partial_{0} M$, we first construct $M(\alpha)$.
The new orbifold, denoted $M(d \alpha)$, is gotten by giving the core curve, $S^{1} \times \{ 0 \}$, in the filling solid torus a transverse angle of $2 \pi / d$, making it (part of) the singular locus with local group ${\mathbb Z} / d {\mathbb Z}$.
We say that $M(d \alpha)$ is obtained from $M$ by {\em $d \alpha$-orbifold Dehn filling}, or simply {\em $d \alpha$-Dehn filling}.
For convenience, we will refer to the positive integer and slope together, $d \alpha$, as a {\em generalized slope}.
As above, if $\mu, \lambda$ is a basis for $H_{1}(\partial_{0}M ; {\mathbb Z})$, and $\alpha$ is given by $(p,q)$, we denote $M(d \alpha)$ by $M(dp,dq)$.

It will often be the case that we wish to fill several boundary components of a compact manifold $M$.
If we have the toroidal boundary components of $M$ labelled $\partial_{1}M,...,\partial_{k}M$, if $d_{j} \alpha_{j}$ is a generalized slope on $\partial_{j} M$ for $j =1,..,k$, then {\em $(d_{1} \alpha_{1},...,d_{k} \alpha_{k})$-Dehn filling on $M$}, denoted $M(d_{1} \alpha_{1},...,d_{k} \alpha_{k})$, is the result of $d_{j} \alpha_{j}$-Dehn filling each of $\partial_{j} M $, for each $j = 1,..., k$.
When we wish to fill only some of the toroidal boundary components, we use the $\infty$ symbol in place of the slope information.
Thus, $M(d_{1}\alpha_{1},\infty, d_{3} \alpha_{3},d_{4}\alpha_{4}, \infty)$ is the orbifold for which $\partial_{1} M$, $\partial_{3} M$ and $\partial_{4} M$ have been filled along $d_{1} \alpha_{1}$, $d_{3} \alpha_{3}$, and $d_{4} \alpha_{4}$ respectively, while $\partial_{2} M$ and $\partial_{5} M$ remain un-filled.

Given a compact $3$-manifold containing a tame link $L \subset M$, we let $N(L)$ denote an open tubular (or regular) neighborhood of $L$.
The {\em exterior of $L$ in $M$} is given by $X_{M}(L) = M \setminus N(L)$.
$\partial X_{M}(L) \setminus \partial M$ is a disjoint union of tori.
In the special case that $M = S^{3}$ we will write $X_{S^{3}}(L) = X(L)$.

\begin{convention} \label{convent}
As is commonly done, we will often make no distinction between the complement of $L$ and the exterior of $L$, denoting both by $X_{M}(L)$.
Quite often this distinction is unimportant.
However, when it is, it should be clear from the context which we are referring to.
This convention greatly simplifies the notation and consequently the exposition.
\end{convention}

\pagebreak

\subsection{Riemannian geometry} \label{Riemanniansection}

We review some terminology and facts from Riemannian geometry.
See \cite{DoCarmo} for more details.\\

Let $M$ be a smooth manifold with a Riemannian metric $g$, Levi-Civita connection $\nabla$, and associated covariant derivative $\frac{D}{dt}$.
Given a point $p \in M$, and a 2-dimensional subspace $\sigma_{p} \subset T_{p}(M)$, we denote the sectional curvature at $\sigma_{p}$ by $K(\sigma_{p})$.
When we wish to emphasize the metric, we may write $K_{g}(\sigma_{p})$.

For any unit speed path $\gamma : [a,b] \rightarrow M$ the {\em geodesic curvature} of $\gamma$ is defined to be the function $\kappa_{\gamma}: [a,b] \rightarrow {\mathbb R}$ given by
$$\kappa_{\gamma} = \sqrt{g \left( \frac{D \gamma}{dt},\frac{D \gamma}{dt} \right) }$$
This function measures the deviation of $\gamma$ from being geodesic:  $\gamma$ is a geodesic if and only if $\kappa_{\gamma} = 0$.

Suppose now that $M$ is 3-dimensional and let $F \subset M$ be an oriented surface in $M$.
The orientations on $F$ and $M$ determined a unique unit normal field $\eta :F \rightarrow TF^{\perp}$.
Here $TF^{\perp}$ is the orthogonal complement to $TF$ in $TM$.

The surface $F$ naturally inherits a Riemannian metric $g_{F}$ for which the inclusion is an isometric embedding.
We denote the Levi-Civita connection on $F$ determined by $g_{F}$ by $\nabla^{F}$, which is given by projecting $\nabla$ onto $TF$.
More precisely, for any vector fields $X, Y$ on $F$, taking any smooth extensions to a neighborhood in $M$, $\overline{X},\overline{Y}$, we have
$$\nabla^{F}_{X}Y = \pi_{T} \nabla_{\overline{X}} \overline{Y}$$
where $\pi_{T} : TM|_{F}  \rightarrow TF$ is the orthogonal projection.

Projecting $\nabla$ onto $TF$ describes the intrinsic geometry of $F$.
By projecting $\nabla$ onto $TF^{\perp}$ we can describe the extrinsic geometry of $F \subset M$.
This is most conveniently measured by the {\em second fundamental form}, $\Pi$, which is a symmetric bilinear form on $TF$ defined by
$$\Pi(X,Y) = g(\nabla_{\overline{X}} \overline{Y},\eta)$$
where $X, Y$ and $\overline{X},\overline{Y}$ are as above.

Together, $\Pi$ and $g_{F}$ dually determine the {\em shape operator}
$$S_{p}: T_{p}F \rightarrow T_{p}F$$
at every $p \in F$, by the formula $g_{F}(S_{p}(X),Y) = \Pi(X,Y)$.
$S_{p}$ is symmetric, and the pair of real eigenvalues of $S_{p}$, $\lambda_{1}(p), \lambda_{2}(p)$, are the {\em principal curvature} of $F$ in $M$ at $p$.
$F$ is {\em totally geodesic} when $\lambda_{1}$ and $\lambda_{2}$ vanish.

If $\gamma : [a,b] \rightarrow F$ is any unit speed geodesic in $F$, then we may view $\gamma$ as a unit speed path into $M$, and as such we can consider its geodesic curvature, $\kappa_{\gamma}$.
This is bounded by
\begin{equation} \label{geodcurvboundeqn}
\kappa_{\gamma}(t) \leq \max \{| \lambda_{1}(\gamma(t))|,| \lambda_{2}(\gamma(t))| \}
\end{equation}

The sectional curvature $K_{g}(T_{p}F)$ and the sectional (or Gaussian) curvature $K_{g_{F}}(T_{p}F)$ are related by the following

\begin{theorem}[Gauss] \label{gauss} 
Suppose $F \subset M$ is as above.
Then for every $p \in F$
$$K_{g_{F}}(T_{p}F) - K_{g}(T_{p}F) = \lambda_{1}(p) \lambda_{2}(p)$$
\end{theorem}

\subsection{Hyperbolic geometry} \label{hypergeomsection}

We recall a few facts from hyperbolic geometry.
See \cite{Tnotes}, \cite{BenPet}, and \cite{Ratcliffe} for more details.\\

{\em Hyperbolic $n$-space}, denoted ${\mathbb H}^{n}$, is the unique (up to isometry) complete, simply connected $n$-manifold with constant sectional curvature -1.
One model for hyperbolic $n$-space, is the {\em upper half space}
$$\{ (x_{1},...,x_{n}) \in {\mathbb R}^{n} \, : \, x_{n} > 0 \}$$
equipped with the metric
$$ds^{2} = \frac{dx_{1}^{2} + ... + dx_{n}^{2}}{x_{n}^{2}}$$
We will make no distinction between hyperbolic space and the upper half space model of hyperbolic space.
The group of isometries of ${\mathbb H}^{n}$ will be denoted by $\Isom({\mathbb H}^{n})$

A {\em hyperbolic $n$-manifold} is a complete Riemannian $n$-manifold with constant sectional curvature -1.
The universal cover of any hyperbolic $n$-manifold is isometric to ${\mathbb H}^{n}$ with the pull-back metric.
Consequently, the covering group of $M$ acts by isometries, and so we may view
$$M = {\mathbb H}^{n} / \Gamma$$
where $\Gamma < \Isom({\mathbb H}^{n})$ is a discrete torsion free group isomorphic to $\pi_{1}(M)$.
We will also consider such quotients in which $\Gamma$ is allowed to have torsion.
In this case the quotient ${\mathbb H}^{n} / \Gamma$ is a {\em hyperbolic $n$-orbifold}.

A {\em horoball} is the image of the set
$$H_{0}  = \{ (x_{1},...,x_{n}) \in {\mathbb H}^{n} \, : \, x_{n} \geq 1 \}$$
under an isometry $\phi \in \Isom({\mathbb H}^{n})$.
The boundary of a horoball is flat with the induced metric and the stabilizer of a horoball in $\Isom({\mathbb H}^{n})$ is isomorphic to the isometry group of Euclidean $(n-1)$-space.

For our purposes, a rank-$k$ horoball cusp (for $k \leq n -1$) is a quotient of a horoball $H \subset {\mathbb H}^{n}$ by a discrete rank $k$ free Abelian subgroup $\Gamma_{H} < Stab_{\Isom({\mathbb H}^{n})}(H)$.

We will be primarily concerned with hyperbolic 2- and 3-manifolds.
A consequence of the Margulis Lemma is that an orientable, finite volume hyperbolic 3-manifold, $M$, is the interior of a compact 3-manifold $\overline{M}$ with (possibly empty) toroidal boundary.
There is a product neighborhood of the boundary of $\overline{M}$ whose intersection with $M$ consists of pairwise disjoint embedded rank-$2$ horoball cusps.
The complement of the interior of these horoball cusp neighborhoods in $M$ is homeomorphic to $\overline{M}$.
We will not make a distinction between $M$ and $\overline{M}$ (see also Convention \ref{convent}).
For example, we may refer to a compact 3-manifold as being hyperbolic, by which we mean that the interior is hyperbolic.
We will also refer to Dehn filling a cusp of $M$, by which we mean Dehn filling the corresponding boundary component of $\overline{M}$.
By a (generalized) slope on a cusp of $M$, we mean a (generalized) slope on the corresponding boundary component of $\overline{M}$

Thurston has shown that ``most'' 3-manifolds are hyperbolic (see \cite{Tbull} and \cite{Otal}).
One instance of this is the following (see \cite{Tnotes} and \cite{HodgKerck})

\begin{theorem}[Thurston] \label{hyperbolicdehnsurgery}
If $M$ is a hyperbolic $3$-manifold of finite volume, with cusps $C_{1},...,C_{k}$, then for each $i = 1,...,k$, there is a finite set, $E_{i}$, of generalized slopes on $C_{i}$, such that orbifold $M(d_{1} \alpha_{1},...,d_{k} \alpha_{k})$ is hyperbolic, provided $d_{i} \alpha_{i} \not \in E_{i}$ for every $i = 1,..., k$.
\end{theorem}

The hyperbolic structures of $M$ and $M(d_{1} \alpha_{1},...,d_{k} \alpha_{k})$ are related.
In particular, an infinite sequence of distinct (on each cusp) Dehn fillings will {\em converge geometrically} to the original manifold $M$.
The following description of this geometric convergence will be necessary for the construction of the examples in Section \ref{knotsection} (see \cite{BenPet}).

For each $i = 1,...,k$, we let $\{ d_{ij}\alpha_{ij} \}_{j=1}^{\infty}$ be an infinite sequence of distinct generalized slopes on the cusp $C_{i}$ (outside the set $E_{i}$).
We write $\pi_{M}:{\mathbb H}^{3} \rightarrow M$ and $\pi_{M_{j}}:{\mathbb H}^{3} \rightarrow M_{j} = M(d_{1j}\alpha_{1j},...,d_{kj}\alpha_{kj})$ to denote the universal covers.

\begin{theorem} \label{geometricconvergence}
There exists smooth embeddings $\phi_{j}:M \rightarrow M_{j}$ and lifts $\widetilde{\phi}_{j} : {\mathbb H}^{3} \rightarrow {\mathbb H}^{3}$, such that for any $q \in {\mathbb H}^{3}$ and $R > 0$, the sequence $\{ {\widetilde{\phi}_{j}}|_{B(q,R)} \}_{j=1}^{\infty}$ converges in the $C^{\infty}$ topology on $B(q,R)$ to the identity.
\end{theorem}

Roughly speaking, this theorem says that as $j \rightarrow \infty$, larger and larger compact subsets of $M$ look more and more like larger and larger compact subsets of $M_{j}$.\\

One proof of Theorem \ref{hyperbolicdehnsurgery}, based on ideal triangulations, is given in \cite{PetPort} (see also \cite{BenPet}).
J. Weeks has written a computer program, {\em SnapPea}, based on ideal triangulations which computes approximate hyperbolic structures on link complements (see \cite{Weeks}).
Although Weeks' program does not provide a rigorous proof of hyperbolicity, recent work of H. Moser \cite{Moser} is able to bridge the computational imprecision in Weeks' program with a quantitative version of the Inverse Function Theorem.
In particular, Weeks' program in conjunction with Moser's (along with O. Goodman's application {\em Snap}) can be used to {\em prove} the hyperbolicity of certain manifolds.
We have appealed to these programs in Sections \ref{linksection} and \ref{knotsection} to find hyperbolic structures, and thank Moser for her time and effort in carrying out those calculations (see \cite{Moser}).

\subsection{Quasi-isometry} \label{quasiisomsection}

There is another well known fact about hyperbolic space which we will need (see \cite{Tnotes} and \cite{Bonahon} for slight variations on this statement).
We have included a proof of this in an appendix at the end of the paper, as a convenience for the reader.

Given numbers $k \geq 1$ and $c \geq 0$, a map $\phi:(X,d) \rightarrow (Y,\rho)$ between metric spaces is a {\em $(k,c)$-quasi-isometry} (or {\em $k$-quasi-isometry} if $c = 0$), if
$$\frac{1}{k}d(x,x') - c \leq d(\phi(x),\phi(x')) \leq k d(x,x') + c$$
If $X$ is a metric space and $I$ is any interval in ${\mathbb R}$, then a $k$-quasi-isometry $\gamma:I \rightarrow X$ is called a {\em $k$-quasi-geodesic}.

\begin{lemma} \label{quasicurvatureclose}
There exists a continuous non-negative function $f$ on the interval $[0,1)$ with $f(0)=0$ and having the following property.
Suppose $\gamma : [a,b] \rightarrow {\mathbb H}^{n}$ is a unit speed path whose geodesic curvature satisfies
$$\kappa_{\gamma}(t) \leq {\cal K}$$
for all $t \in [a,b]$, where $0 \leq {\cal K} < 1$ is some constant.
Then $\gamma$ is a $\frac{1}{\sqrt{1-{\cal K}^{2}}}$-quasi-geodesic.
Moreover, if $g_{\gamma}:[a,b] \rightarrow {\mathbb H}^{n}$ is the unique geodesic connecting the endpoints $\gamma(a)$ and $\gamma(b)$, then the Hausdorff distance between the image of $\gamma$ and $g_{\gamma}$ is no more than $f({\cal K})$.
\end{lemma}

\subsection{Surfaces and hyperbolic geometry}

A closed surface $F$ in a hyperbolic manifold $M$ is {\em quasi-Fuchsian} if a lift of the inclusion $\widetilde{F} \rightarrow \widetilde{M}$ of universal covers is a quasi-isometry.
Here we are using the pull back metric on $\widetilde{F}$ (with respect to any metric on $F$).
This is easily seen to be equivalent to the usual notion of quasi-Fuchsian for closed surfaces.\\

To guarantee that a surface in a hyperbolic 3-manifold is totally geodesic, we often use the following (see \cite{MenascoReid})

\begin{lemma} \label{fixedpointset}
If $M$ is a finite volume hyperbolic 3-orbifold and $\phi : M \rightarrow M$ is an orientation reversing involution fixing a $2$-orbifold $F$, then $F$ can be homotoped to be totally geodesic.
\end{lemma}

\noindent
{\em Sketch of proof.} It follows from the irreducibility of $M$ that $F$ is incompressible.
We can find a lift of $\phi$ to the universal cover
$$\widetilde{\phi} : {\mathbb H}^{3} \rightarrow {\mathbb H}^{3}$$
which fixes a component, $\widetilde{F}$, of the preimage of $F$ pointwise.
The map $\widetilde{\phi}$ extends continuously to the sphere at infinity
$$\partial\widetilde{\phi}: S^{2}_{\infty} \rightarrow S^{2}_{\infty}$$
As in the proof of the Mostow Rigidity Theorem (see e.g. \cite{Tnotes}), it follows that $\partial \widetilde{\phi}$ is the extension of an isometry.
Moreover, since $\phi^{2} = id_{M}$ and by our choice of lift $\widetilde{\phi}$, we see that $\widetilde{\phi}^{2} = id_{{\mathbb H}^{3}}$.
Therefore $\partial \widetilde{\phi}^{2} = id_{S^{2}_{\infty}}$.

It follows that $\partial \widetilde{\phi}$ fixes a geometric circle which must be the boundary of $\widetilde{F}$ at infinity.
Therefore, $\pi_{1}(F) \subset \pi_{1}(M)$ (acting by covering transformations) must stabilize this circle and hence $F$ can be homotoped to be totally geodesic. \hfill $\Box$\\

Another source of totally geodesic surfaces come from {\em triangle orbifolds}.
A triangle orbifold is a 2-orbifold which is topologically an $n$-times punctured sphere, for $n = 0,1,2,3$, with $3-n$ cone points.
A $(p_{1},p_{2},p_{3})$-triangle orbifold is a triangle orbifold where the cone points have orders given by $p_{i}$, for $i=1,2,3$ (and $p_{i} = \infty$ means that the instead of a cone point, one has a puncture).
For example, a $(2,3,\infty)$-triangle orbifold is a once punctured sphere with one cone point of order $2$ and one of order $3$.
Any triangle suborbifold of a hyperbolic 3-orbifold is incompressible and moreover, by applying an isotopy, we may assume it is totally geodesic (this follows as in \cite{Adams}).
These will also be useful in applying the cut-and-paste techniques of \cite{Adams}.

\section{Links} \label{linksection}

In this section we prove the following

\begin{theorem} \label{links}
For any even integer $g \geq 2$, there exists a two component hyperbolic link in $S^{3}$ which contains an embedded totally geodesic surface of genus $g$ in its complement.
\end{theorem}

\begin{figure}[htb]
\begin{center}
\ \psfig{file=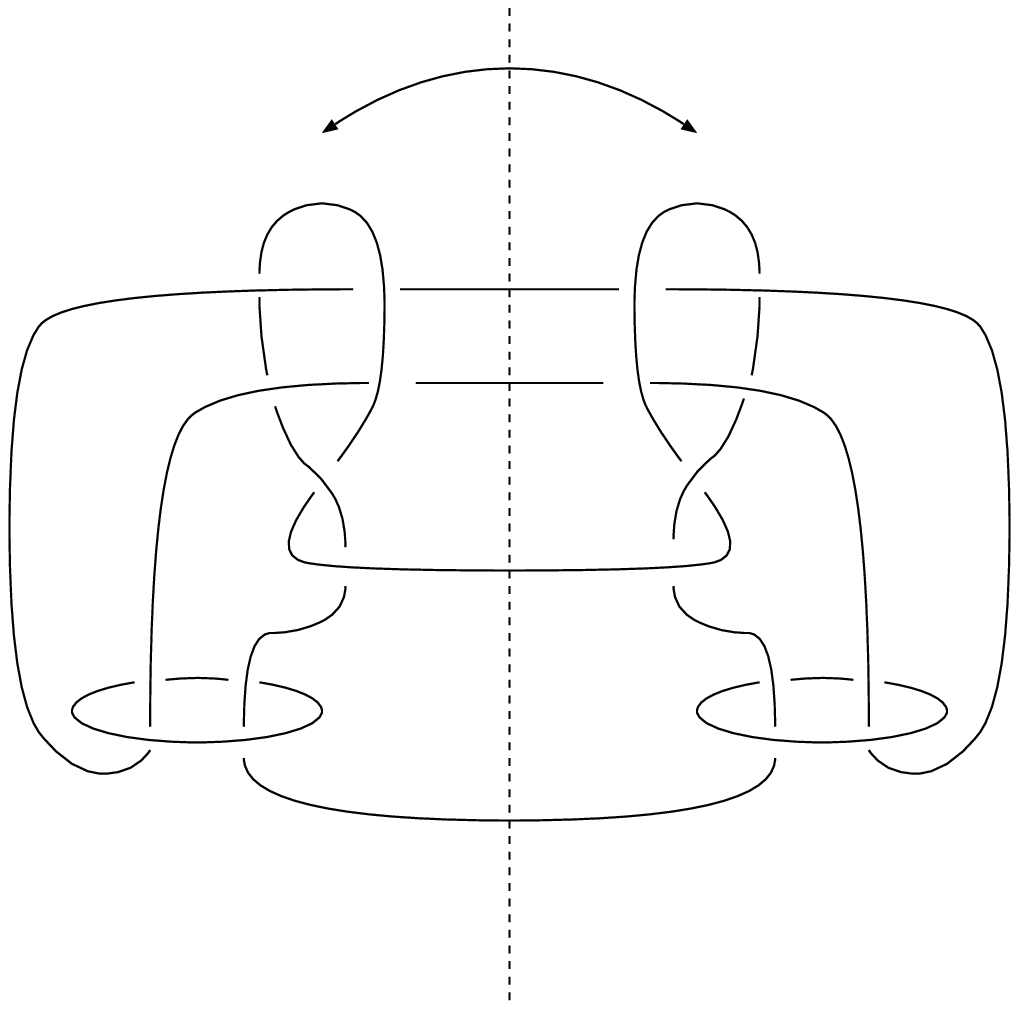,height=1.9truein}
\caption{$L_{0}$ and the involution $\phi$.}
\label{linksurface1}
\end{center}
  \setlength{\unitlength}{1in}
  \begin{picture}(0,0)(0,0)
    \put(2.7,2.4){$\phi$}
    \put(1.9,1.9){$K_{3}$}
    \put(3.5,2.1){$K_{4}$}
    \put(2.68,1.15){$K_{1}$}
    \put(3.18,1.15){$K_{2}$}
  \end{picture}
\end{figure}

\noindent
{\em Proof. } We begin with the link $L_{0}$ shown in Figure \ref{linksurface1}, with components labelled $K_{1},...,K_{4}$ as indicated.
According to {\em SnapPea}, and verified by Moser \cite{Moser} (see \S \ref{hypergeomsection}), $M_{0}=X(L_{0})$ admits a hyperbolic structure.
Further, $M_{0}$ admits an orientation reversing involution $\phi : M_{0} \rightarrow M_{0}$, as indicated in Figure \ref{linksurface1}, fixing a 4-punctured sphere which is thus totally geodesic by Lemma \ref{fixedpointset}.

By Theorem \ref{hyperbolicdehnsurgery}, $M_{0}(\infty,\infty,(p,0),(p,0))$ is a hyperbolic orbifold for any sufficiently large positive integer $p$.
In fact, appealing to the geometrization theorem for orbifolds (see \cite{CKH} and \cite{BLP}) and some topology, this holds for all $p \geq 3$.
We first check that $M_{0}(\infty, \infty, (p,0), (p,0))$ is orbifold irreducible with orbifold incompressible boundary and contains no essential euclidean $2$-orbifold.
All of these amount to showing that there are no $2$-orbifolds with certain properties.
This is easy to check since we must only consider honest $2$-orbifolds; an offending surface would live in the complement of the singular locus and would so give rise to an offending surface after drilling out the singular locus, contradicting the hyperbolicity of $M_{0}$.
It follows that for $p \geq 3$, $M_{0}(\infty,\infty,(p,0),(p,0))$ has a geometric structure which must be hyperbolic (one can easily rule out any Seifert fibered structure).

The involution $\phi$ persists in the filled manifold, and so $M_{0}(\infty,\infty,(p,0),(p,0))$ contains a totally geodesic 2-orbifold, $F$, which is a 2-sphere with four order $p$ cone points, again by Lemma \ref{fixedpointset}.

\begin{figure}[htb]
\begin{center}
\ \psfig{file=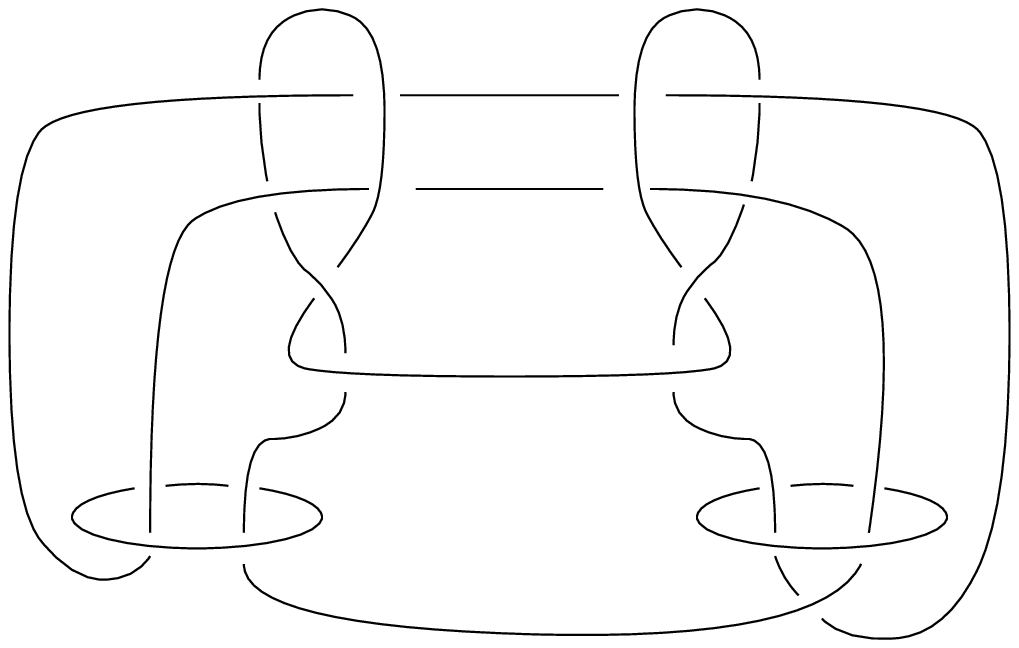,height=1.3truein}
\caption{The link $L_{1}$.}
\label{linksurface3}
\end{center}
  \setlength{\unitlength}{1in}
  \begin{picture}(0,0)(0,0)
    \put(1.8,1.65){$K_{3}$}
    \put(2.65,.85){$K_{1}$}
    \put(3.2,.85){$K_{2}$}
  \end{picture}
\end{figure}

Next, consider the link $L_{1}$ shown in Figure \ref{linksurface3}, with components labelled $K_{1}',K_{2}',K_{3}'$.
$M_{1} = X(L_{1})$ is also hyperbolic as it is obtained from $M_{0}$ by cutting open along a thrice-punctured sphere and gluing back with a half twist, (see \cite{Adams}).
In fact, there is an isometry from the complement of the thrice-punctured sphere in $M_{0}$ to the complement of the thrice-punctured sphere in $M_{1}$.

Similarly, $M_{1}(\infty,\infty,(p,0))$ can be obtained from $M_{0}(\infty,\infty,(p,0),(p,0))$ by cutting open along a $(p,p,\infty)$-triangle orbifold, $T$, and gluing back with a half-twist.
Again, there is an isometry from the complement of $T$ in $M_{0}(\infty,\infty,(p,0),(p,0))$ to the complement of $T$ in $M_{1}(\infty,\infty,(p,0))$.
Since $F$ is disjoint from $T$, its image in $M_{1}(\infty,\infty,(p,0))$ is totally geodesic.

We let $M_{p}$ denote the $p$-fold cyclic branched cover of $S^{3}$, branched over $K_{3}'$, with the preimage of $K_{1}'$ and $K_{2}'$ deleted.
$M_{p}$ is a manifold cover of $M_{1}(\infty,\infty,(p,0))$, and since $K_{3}'$ is unknotted, $M_{p}$ is a link complement in $S^{3}$.
Moreover, since $\lk(K_{1}',K_{3}')=2=\lk(K_{2}',K_{3}')$, it follows that for positive odd integers $p$, the preimages of $K_{1}'$ and $K_{2}'$ in the branched cover are connected.
That is, $M_{p} \cong X(L_{p})$, where $L_{p} \subset S^{3}$ is a two-component link.
The preimage of $F$ in $M_{p}$ is thus a closed totally geodesic surface in a two-component link of genus $p-1$. \hfill $\Box$\\

We have drawn the link $L_{7}$ in Figure \ref{holycrap}.

\section{Knots} \label{knotsection}

In this section we construct a family of knots proving the next two theorems.

\begin{theorem} \label{knots}
For any $g \geq 3$ and any $\epsilon > 0$, there exists a hyperbolic knot $K \subset S^{3}$ containing a closed embedded surface of genus $g$ in its complement whose principal curvatures are bounded in absolute value by $\epsilon$.
\end{theorem}

\begin{theorem} \label{pinchknots}
For any $g \geq 3$ and any $\epsilon >0$, there exists a knot $K \subset S^{3}$ with complement supporting a Riemannian metric with negative sectional curvatures pinched between $-1-\epsilon$ and $-1+\epsilon$, which contains a closed embedded totally geodesic surface of genus $g$.
\end{theorem}

These two theorems are both consequences of the following lemma.

\begin{lemma} \label{knotlemma}
For any $g \geq 3$, there exists a sequence of knots $\{ K_{j} \}_{j=1}^{\infty}$ and a link $L$ such that
\begin{enumerate}
\item $X(L)$ and each $X(K_{j})$ are hyperbolic,
\item $X(K_{j}) \cong X(L)((p_{1,j},q_{1,j}),...,(p_{k,j},q_{k,j}),\infty)$, for an infinite sequence of distinct slopes, $\{p_{i,j},q_{i,j}\}_{j=1}^{\infty}$, on the $i^{th}$ cusp of $X(L)$,
\item $X(L)$ contains a closed totally geodesic surface $F$ with genus $g$.
\end{enumerate}
\end{lemma}

The construction which proves this lemma is deferred to \S \ref{knotsandlinksection}.
In \S \ref{knotsection1} and \S \ref{knotsection2} we use Lemma \ref{knotlemma} to prove Theorems \ref{knots} and \ref{pinchknots}, respectively.

\subsection{Small curvature surfaces in hyperbolic knot complements} \label{knotsection1}

\noindent
{\em Proof of Theorem \ref{knots}.} We consider the surface $F$, the family of knots $\{K_{j}\}_{j=1}^{\infty}$, and the link $L$ from Lemma \ref{knotlemma}.
By Theorem \ref{geometricconvergence}, we have embeddings $\phi_{j}:X(L) \rightarrow X(K_{j})$ and lifts $\widetilde{\phi}_{j}:{\mathbb H}^{3} \rightarrow {\mathbb H}^{3}$ satisfying Theorem \ref{geometricconvergence}.
The embeddings $\phi_{j}$ restrict to embeddings
$${\phi_{j}}|_{F} : F \rightarrow X(K_{j})$$

We will be done if we can show that the principal curvatures of ${\phi_{j}}|_{F}$ converge to $0$ as $j \rightarrow \infty$.
We use the notation of Theorem \ref{geometricconvergence}, replacing $M$ by $X(L)$ and $M_{j}$ by $X(K_{j})$.\\

Since $F$ is totally geodesic, there is a totally geodesic hyperbolic plane ${\mathbb H}^{2} \subset {\mathbb H}^{3}$ covering $F$.
The restriction of the universal cover of $X(L)$
$${\pi_{X(L)}}|_{{\mathbb H}^{2}} : {\mathbb H}^{2} \rightarrow F$$
is the universal cover of $F$.
Since $F$ is compact, its diameter is finite.
Therefore, there exists $R > 0$ and $q \in {\mathbb H}^{2}$, such that $B(q,R) \cap {\mathbb H}^{2}$ contains a fundamental domain for the action of $\pi_{1}(F)$.
In particular, $\pi_{X_{L}}(B(q,R) \cap {\mathbb H}^{2}) = F$.

Now note that
$${\widetilde{\phi}_{j}}|_{B(q,R) \cap {\mathbb H}^{2}}$$
is converging in the $C^{\infty}$ topology to a totally geodesic embedding.
Therefore, the second fundamental forms $\widetilde{\Pi}_{j}$ for ${\widetilde{\phi}_{j}}|_{B(q,R) \cap {\mathbb H}^{2}}$ are converging to zero uniformly on $B(q,R) \cap {\mathbb H}^{2}$.

Moreover, since
$$\pi_{X(K_{j})} \circ {\widetilde{\phi}_{j}}|_{B(q,R) \cap {\mathbb H}^{2}} = {\phi_{j}}|_{F} \circ {\pi_{X(L)}}|_{B(q,R) \cap {\mathbb H}^{2}}$$
and $\pi_{X(L)}$ and $\pi_{X(K_{j})}$ are local isometries, it follows that the second fundamental forms $\Pi_{j}$ for ${\phi_{j}}|_{F}$ are converging uniformly to zero on $F$.
That is, the principal curvatures of the embeddings of $F$ into $X(K_{j})$ are converging to zero as required. \hfill $\Box$\\

\subsection{Totally geodesic surfaces in nearly hyperbolic knot complements} \label{knotsection2}

\noindent
{\em Proof of Theorem \ref{pinchknots}.} Again, we let $F$, $\{ K_{j} \}_{j=1}^{\infty}$, and $L$ be as in Lemma \ref{knotlemma}.

If we let $C_{1},..., C_{k}$ denote embedded horoball cusps in $X(L)$ which are pairwise disjoint and also disjoint from $F$, then as we let $j \rightarrow \infty$, the lengths of the filling curves on the boundary of $C_{1} \cup ... \cup C_{k}$ are approaching infinity.

As in the proof of the Gromov-Thurston $2\pi$-Theorem (see \cite{BleiHodg}), for each of the fillings we seek to construct Riemannian metrics on solid tori $V_{1},...,V_{k}$ with negative sectional curvature for which a neighborhood of the boundary of $V_{i}$ is isometric to a neighborhood of the boundary of the cusp $C_{i}$ by an isometry taking the boundary of the meridian of $V_{i}$ to the filling curve of $C_{i}$, for each $i=1,...,k$.
Moreover, we wish to do this so that given $\epsilon > 0$, for all sufficiently large $j$, the metric on each of the $V_{j}$'s is pinched between $-1-\epsilon$ and $-1+\epsilon$.
We may then perform the Dehn fillings (with $j$ sufficiently large) requiring that the gluing of $V_{1},...,V_{k}$ to $X(L) \setminus int(C_{1} \cup \cdots \cup C_{k})$ to be by isometries.
The resulting manifold, which is diffeomorphic to $X(K_{j})$, is equipped with a Riemannian metric having sectional curvatures pinched between $-1-\epsilon$ and $-1+\epsilon$, and moreover, the metric on
$$X(L) \setminus \bigcup_{i=1}^{k} C_{j} \subset X(K_{j})$$
has not changed, and so still contains the closed embedded totally geodesic surface $F$.

Thus, to complete the proof of Theorem \ref{pinchknots}, we need only prove the following lemma, analogous to Lemma 10 of \cite{BleiHodg}.

\begin{lemma} \label{metricsolidtori}
There exists a constant ${\cal L} > 0$, so that if $C$ is a rank-2 horoball cusp and $\gamma$ is a geodesic on $\partial C$ with length $l \geq e^{3} \pi$ then there is a metric on a solid torus $V$ such that the $1$-neighborhood of $\partial V$ is isometric to the $1$-neighborhood of $\partial C$ by an isometry taking the boundary of some meridian to $\gamma$.
Furthermore, the sectional curvatures $K(\sigma_{p})$ of $V$ satisfy
$$-1 - \frac{{\cal L}}{l^{2}} \leq K(\sigma_{p}) \leq -1 + \frac{{\cal L}}{l^{2}}$$
for all $p \in V$, $\sigma_{p} \subset T_{p}(V)$.
\end{lemma}

\noindent
{\em Proof. } Rather than constructing our metric on $V$, we will construct a metric on its universal cover
$$\widetilde{V} \cong D^{2} \times {\mathbb R}$$
so that it is invariant under the obvious $S^{1} \times {\mathbb R}$ action.
Further, we will require that in the $1$-neighborhood of the boundary, the metric is isometric to the $1$-neighborhood of a rank-$1$ horoball cusp, and the meridian, $\partial D^{2} \times \{ * \}$, has length $l$.
The lemma will follow then by taking the quotient of $\widetilde{V}$ by an appropriate isometric ${\mathbb Z}$ action.

We use the notation of the proof of Lemma 10 of \cite{BleiHodg} and so consider a metric of the form
$$ds^{2} = dr^{2} + f^{2}(r) d\mu^{2} + g^{2}(r) d\lambda^{2}$$
on $\widetilde{V}$ in cylindrical coordinates $r,\mu,\lambda$; where $r \leq 0$, is the (signed) radial distance measured outwards from $\partial \widetilde{V}$ (so points on $int(\widetilde{V})$ have a negative $r$ coordinate), $0 \leq \mu \leq 1$ is measured in the meridional direction, and $-\infty < \lambda < \infty$ is measured in the direction perpendicular to $\mu$ and $r$.

We recall the following facts from \cite{BleiHodg}
\begin{itemize}
\item if $f$ and $g$ satisfy $f(r) = le^{r}$ and $g(r) = e^{r}$ for $-\epsilon \leq r \leq 0$, then the metric in an $\epsilon$-neighborhood of the boundary is isometric to the $\epsilon$-neighborhood of the boundary of a rank-$1$ horoball cusp, and the meridian curve has length $l$.
\item if the core occurs at $r = r_{0}$, and $f(r) = 2 \pi \sinh(r-r_{0})$ and $g(r) = b \cosh(r-r_{0})$ for $r_{0} \leq r \leq r_{0} + \epsilon$ and for some constant $b > 0$, then the $\epsilon$-neighborhood of the core is non-singular.
\item the sectional curvatures of $ds^{2}$ are convex combinations of the functions
$$\left\{ -\frac{f''}{f} \, \, , \, \, -\frac{g''}{g} \, \, , \, \, -\frac{f' \cdot g'}{f \cdot g} \right\}$$
\end{itemize}

Now, let $\phi : {\mathbb R} \rightarrow {\mathbb R}$ be a smooth function with $0 \leq \phi(r) \leq 1$ satisfying
$$\phi(r) = \left\{ \begin{array}{ll}
1 & \mbox{ for } r \leq -2\\
0 & \mbox{ for } r \geq -1\\ \end{array} \right.$$

Set $r_{0} = -\log \left( \frac{l}{\pi} \right)$ and note that $r_{0} \leq -3$ since $l \geq e^{3} \pi$.
For $r_{0} \leq r \leq 0$, define
$$f(r) = \pi(e^{r-r_{0}} - \phi(r) e^{r_{0} - r}) \mbox{ and } g(r) = e^{r} + \phi(r) e^{2r_{0} - r}$$
and note that for $ -1 \leq r \leq 0$, we have
$$f(r) = \pi e^{-r_{0}} e^{r} = l e^{r} \mbox{ and } g(r) = e^{r}$$
and for $r_{0} \leq r \leq -2$, we have
$$f(r) = 2 \pi \frac{e^{r - r_{0}} - e^{r_{0} - r}}{2} = 2 \pi \sinh(r - r_{0})$$
and
$$g(r) = 2e^{r_{0}} \left( \frac{e^{r - r_{0}} + e^{r_{0} - r}}{2} \right) = 2e^{r_{0}} \cosh(r-r_{0})$$

From the calculations of \cite{BleiHodg} mentioned above, we will be done if we can show that
\begin{equation} \label{sectionalcurvatures}
|-1 - \left( -\frac{f''}{f} \right) |, |-1 - \left( -\frac{g''}{g} \right) |, |-1 - \left( -\frac{f' \cdot g'}{f \cdot g} \right) | \leq \frac{{\cal L}}{l^{2}}
\end{equation}
for some ${\cal L} > 0$.

By inspection, we see that for $-1 \leq r \leq 0$ and $r_{0} \leq r \leq -2$, we have
$$-\frac{f''(r)}{f(r)} = -1 \, \, , \, \, -\frac{g''(r)}{g(r)} = -1 \, \, , \mbox{ and } -\frac{f'(r) \cdot g'(r)}{f(r) \cdot g(r)} = -1$$
Thus, to verify (\ref{sectionalcurvatures}), we need only check this for $-2 \leq r \leq -1$.

For $-2 \leq r \leq -1$, we have
$$| -1 - \left( -\frac{f''(r)}{f(r)} \right) | = |-1 + \frac{\pi(e^{r-r_{0}} - e^{r_{0} - r}(\phi(r) + \phi''(r) - 2 \phi'(r)))}{\pi(e^{r-r_{0}} - e^{r_{0} - r}\phi(r))} | = \frac{ | e^{r_{0} - r}(\phi''(r) - 2 \phi'(r))|}{|e^{r-r_{0}} - e^{r_{0} -r}\phi(r)|}$$
$$= e^{2r_{0}} e^{-2r} \frac{|\phi''(r) - 2 \phi'(r)|}{| 1 - e^{2(r_{0} - r)} \phi(r)|} =  e^{-2 \log(\frac{l}{\pi})} e^{-2r} \frac{|\phi''(r) - 2 \phi'(r)|}{| 1 - e^{2(r_{0} - r)} \phi(r)|} \leq \frac{\pi^{2}}{l^{2}} e^{4} \frac{|\phi''(r) - 2 \phi'(r)|}{1 - e^{-2}}$$

A similar calculation, shows
$$|-1 - \left( -\frac{g''(r)}{g(r)} \right)| = e^{2(r_{0} - r)} \frac{|\phi''(r) - 2\phi'(r)|}{|1+ \phi(r) e^{2(r_{0} - r)}|} \leq \frac{\pi^{2}}{l^{2}} e^{4} |\phi''(r) - 2 \phi'(r)|$$
and
$$|-1 - \left( -\frac{f'(r) \cdot g'(r)}{f(r) \cdot g(r)} \right) | = e^{4(r_{0} - r)} \frac{|2\phi(r)\phi'(r) - (\phi'(r))^{2}|}{|1 - \phi^{2}(r)e^{4(r_{0} - r)}|} \leq \frac{\pi^{4}}{l^{4}} e^{8} \frac{|2 \phi(r) \phi'(r) - (\phi'(r))^{2}|}{1 - e^{-4}}$$
for $-2 \leq r \leq -1$.

It follows that (\ref{sectionalcurvatures}) holds if we set
$${\cal L} = \max \left\{ \pi^{2} e^{4} \frac{|\phi''(r) - 2 \phi'(r)|}{1 - e^{-2}} \, , \pi^{4} e^{8} \frac{|2 \phi(r) \phi'(r) - (\phi'(r))^{2}|}{1 - e^{-4}} \, \, : \, \, -2 \leq r \leq -1 \, \, \right\}$$
\hfill $\Box$

\subsection{An interesting link} \label{knotsandlinksection}

\begin{figure}[htb]
\begin{center}
\ \psfig{file=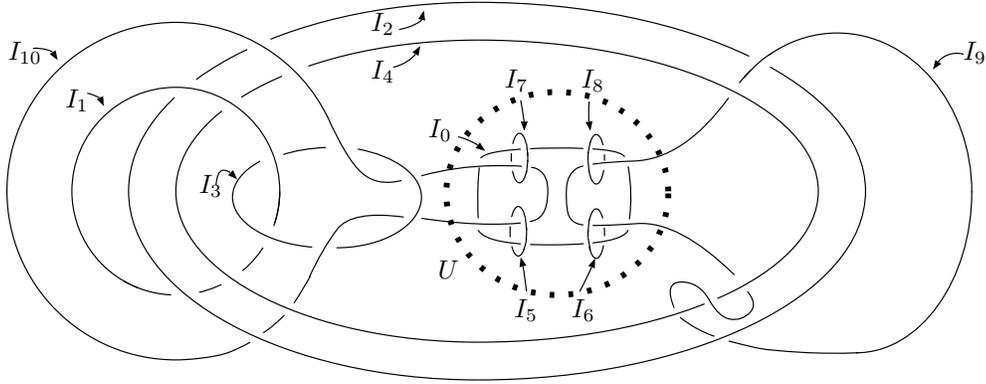,height=2truein}
\caption{The link $J$.}
\label{gettintheknot}
\end{center}
  \setlength{\unitlength}{1in}
  \begin{picture}(0,0)(0,0)
    \put(.5,2.3){$I_{10}$}
    \put(5.5,2.3){$I_{9}$}
    \put(.8,2.05){$I_{1}$}
    \put(1.5,1.6){$I_{3}$}
    \put(2.4,2.45){$I_{2}$}
    \put(2.4,2.2){$I_{4}$}
    \put(3.1,2.15){$I_{7}$}
    \put(3.5,2.15){$I_{8}$}
    \put(3.15,.95){$I_{5}$}
    \put(3.45,.95){$I_{6}$}
    \put(2.7,1.89){$I_{0}$}
    \put(2.75,1.15){$U$}
  \end{picture}
\end{figure}

{\em Proof of Lemma \ref{knotlemma}.} Start with the link $J$ shown in Figure \ref{gettintheknot}, with components $I_{0},...,I_{10}$ as indicated (ignore the dotted circle and the region $U$ bounded by it for the moment).
According to {\em SnapPea}, and verified by Moser \cite{Moser} (see \S \ref{hypergeomsection}), $X(J)$ is hyperbolic.
We refer to the components of $\partial X(J)$ by $\partial_{i} X(J)$ so that $\partial_{i} X(J)$ is the boundary of a neighborhood of $I_{i}$, for $i = 0,...,10$.
For each $i$, let $m_{i},l_{i}$ denote a standard basis for $H_{1}(\partial_{i} X(J))$.

$X(J)$ admits an orientation reversing involution fixing a twice-punctured torus.
To see this, we first note that $I_{1} \cup I_{2}$ is a Hopf link and hence $X(I_{1} \cup I_{2}) \cong T^{2} \times [0,1]$ admits an orientation reversing involution
$$\tau: X(I_{1} \cup I_{2}) \rightarrow X(I_{1} \cup I_{2})$$
fixing a torus, $T \cong T^{2} \times \{ \frac{1}{2} \}$.
We then add components $I_{3}$, $I_{5}$, $I_{7}$, and $I_{9}$ in the complement of $T$ and their respective images $\tau(I_{3}) = I_{4}$, $\tau(I_{5}) = I_{6}$, $\tau(I_{7}) = I_{8}$, and $\tau(I_{9}) = I_{10}$.
Finally, we add the component $I_{0}$ which is invariant under $\tau$ and transversely intersects $T$ twice.
We let $T^{*}$ denote the twice-punctured torus fixed by $\tau$.

Theorem \ref{hyperbolicdehnsurgery} implies that there exists $p_{0} > 0$ such that for any $p > p_{0}$, the orbifold
$${\cal O}_{p} = X(J)((p,0),\infty,...,\infty)$$
is hyperbolic.
In fact, arguing as in \S \ref{linksection}, it suffices to take $p_{0} = 3$.
Moreover, the involution $\tau$ persists in ${\cal O}_{p}$, and so this orbifold contains a totally geodesic 2-orbifold, $T_{p}^{*}$, which is a totally geodesic torus with two cone points of order $p$, by Lemma \ref{fixedpointset}.

The idea for the remainder of the proof is the following.
We find infinitely many distinct fillings on all boundary components of ${\cal O}_{p}$ except $\partial_{10} {\cal O}_{p}$.
For each of these fillings the result will be an orbifold having underlying topological space the complement of a knot $K_{j}$ in $S^{3}$ and singular locus an unknotted curve with cone angle $\frac{2 \pi}{p}$.
Furthermore, the linking number of the singular locus with $K_{j}$ will be one.
Then, the $p$-fold branched cover of $S^{3}$ branched over the singular locus defines a manifold cover of our orbifold which is a knot complement in $S^{3}$.
Infinitely many of these knots can be seen to be obtained by Dehn filling on all but one component of a single link complement (along distinct slopes) which is itself a manifold cover of ${\cal O}_{p}$.
This therefore contains the preimage of $T_{p}^{*}$ which is totally geodesic.\\

Now for $\vec{r} = (r_{1},...,r_{5}) \in {\mathbb Z}^{5}$ consider the orbifold ${\cal O}_{p}(\vec{r})$ defined by
\begin{equation} \label{slopeseqn}
{\cal O}_{p}(\vec{r}) = {\cal O}_{p} ((1+r_{1},-r_{1}),(1,-r_{2}),(1-r_{1},r_{1}),(1,r_{2}),(1,r_{3}),(1,r_{4}),(1,-r_{3}-1),(1,-r_{4}),(1,r_{5}),\infty)
\end{equation}
By Theorem \ref{hyperbolicdehnsurgery} there exists $R > 0$ such that ${\cal O}_{p}(\vec{r})$ is hyperbolic whenever each $|r_{i}| > R$.

For each $i_{1},...,i_{k} \in \{ 0,...,10 \}$, consider the sublink $J_{i_{1},...,i_{k}} \subset J$, given by $$J_{i_{1},...,i_{k}} = J \setminus \left( I_{i_{1}} \cup \cdots \cup I_{i_{k}} \right)$$
So, {\em the indices tell us which components have been left out} (this set is generally smaller than its complement, which is the reason we have chosen this notation).
We will consider the obvious inclusion
$$X(J) \subset M_{i_{1},...,i_{k}}$$
where $M_{i_{1},...,i_{k}} = X(J_{i_{1},...,i_{k}})$.
In particular, we will use this to identify slopes on the boundary components of $M_{i_{1},...,i_{k}}$ with those on the corresponding boundary components of $X(J)$.
For each $\vec{r} \in {\mathbb Z}^{5}$, we let $M_{i_{1},...,i_{k}}(\vec{r})$ denote the manifold obtained from $M_{i_{1},...,i_{k}}$ by filling those boundary components in common with $X(J)$ according to (\ref{slopeseqn}).

We consider $I_{i_{1}} \cup \cdots \cup I_{i_{k}}$ (or any sublink of this) as a link in each of $S^{3}$, $M_{i_{1},...,i_{k}}$, and $M_{i_{i},...,i_{k}}(\vec{r})$.
We sometimes express the dependence on $\vec{r}$ by denoting the component $I_{i}$ in this last manifold by $I_{i}(\vec{r})$.

Note that
\begin{equation} \label{fillingeqn}
X_{M_{0,10}(\vec{r})}(I_{0}(\vec{r}) \cup I_{10}(\vec{r}))((p,0),\infty) \cong {\cal O}_{p}({\vec{r}})
\end{equation}
The slopes are defined in terms of $m_{i},l_{i}$ via the inclusion $X(J) \subset X_{M_{0,10}(\vec{r})}(I_{0} \cup I_{10})$.\\

\noindent
{\bf Claim} For every $\vec{r} \in {\mathbb Z}^{5}$, we have
\begin{enumerate}
\item $M_{0,10}(\vec{r}) \cong S^{3}$.
\item $I_{0}(\vec{r})$ is unknotted.
\item $\lk(I_{0}(\vec{r}),I_{10}(\vec{r})) = 1$
\end{enumerate}

\begin{figure}[htb]
\begin{center}
\ \psfig{file=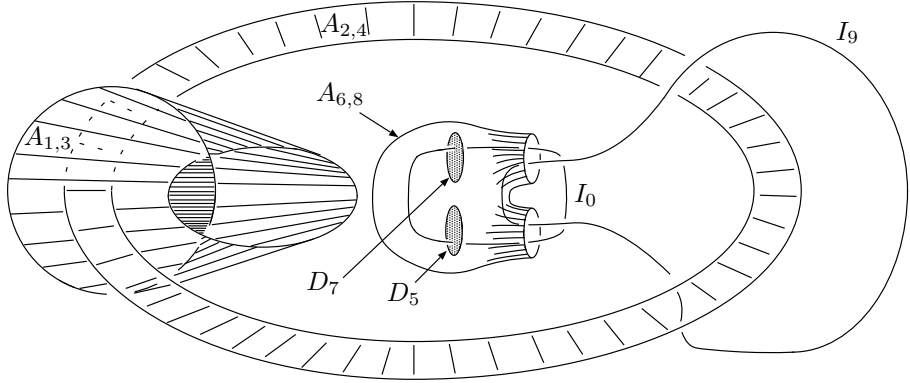,height=2truein}
\caption{annuli and disks with boundaries on $J_{0,9,10}$ ($I_{0}$ and $I_{9}$ are also pictured).}
\label{annuli}
\end{center}
  \setlength{\unitlength}{1in}
  \begin{picture}(0,0)(0,0)
    \put(2.23,1.09){$D_{7}$}
    \put(2.28,2.08){$A_{6,8}$}
    \put(.75,1.85){$A_{1,3}$}
    \put(2.65,1.03){$D_{5}$}
    \put(2.3,2.44){$A_{2,4}$}
    \put(3.63,1.55){$I_{0}$}
    \put(5,2.4){$I_{9}$}
  \end{picture}
\end{figure}

\noindent
{\em Proof of Claim. } We first note that there are annuli $A_{1,3}$, $A_{2,4}$, and $A_{6,8}$ having boundaries $I_{1} \cup I_{3}$, $I_{2} \cup I_{4}$, and $I_{6} \cup I_{8}$, respectively, and disks $D_{5}$ and $D_{7}$ bounded by $I_{5}$ and $I_{7}$, respectively, as shown in Figure \ref{annuli} (the link shown is $J_{10}$).
For any one of these annuli, $A_{i,j}$, we can view it as being embedded in $X(I_{i} \cup I_{j})$ (or $D_{i}$ in $X(I_{i})$).
As such, we can {\em Dehn twist} along $A_{i,j}$ and produce a different embedding of $X(J_{0,9,10})$ into $S^{3}$.
A Dehn twist along an annulus in a three manifold is a homeomorphism of the three manifold supported in a regular neighborhood of the annulus in which one cuts open, twists, and reglues-- if we view the neighborhood of the annulus as $S^{1} \times [0,1] \times [0,1]$, then this is the usual notion of Dehn twist on the $S^{1} \times [0,1]$ factor and the identity on the last $[0,1]$ factor.
In a similar fashion we can Dehn twist along the disks $D_{5}$ and $D_{7}$ to define embeddings of $X(J_{0,9,10})$ into $S^{3}$.

For each such embedding determined by $A_{i,j}$ (respectively, $D_{i}$) we obtain different curves on $\partial_{i}X(J_{0,9,10}) \cup \partial_{j}X(J_{0,9,10})$ (respectively, $\partial_{i}X(J_{0,9,10})$) which bound meridian disks.
In particular, Dehn filling along these curves will result in $S^{3}$ (see \cite{Gordon} for more on this method of altering links while keeping their complements the same).

The new curves bounding meridian disks are described as follows.
If
$$\partial A_{i,j} \cap \partial_{i} X(J_{0,9,10}) = x_{i} \in H_{1}(\partial_{i} X(J_{0,9,10}))$$
then after the $r^{th}$ iterate of the twist in $A_{i,j}$, the new curve bounding a meridian disk on $\partial_{i}X(J_{0,9,10})$ is $m_{i} + r(x_{i} \cdot m_{i})x_{i}$ ("$\cdot$" denote algebraic intersection number).
Likewise, if $\partial A_{i,j} \cap \partial_{j} X(J_{0,9,10}) = x_{j}$, then the new curve bounding a meridian disk on $\partial_{j} X(J_{0,9,10})$ is $m_{j} - r(x_{j} \cdot m_{j})x_{j}$ (we have made an arbitrary choice of direction in which to twist and orientation on the annuli).
In the case of $D_{i}$, the new curve bounding a meridian disk on $\partial_{i}X(J_{0,9,10})$ after the $r^{th}$ iterate of the twist is $m_{i} + rl_{i}$.

Computing the boundary slopes of the annuli, and letting $r,r_{1},...,r_{4} \in {\mathbb Z}$ be any integers (which tell us how many times to twist), we see that
$$X(J_{0,9,10})((1+r_{1},-r_{1}),(1,-r_{2}),(1-r_{1},r_{1}),(1,r_{2}),(1,r_{3}),(1,r_{4}),(1,r),(1,-r_{4}))$$
is $S^{3}$.
We will additionally require that $r = -r_{3}-1$.
Although this is not necessary at the moment, we make this assumption now.

There is one point where we must be a little careful.
Note that $A_{1,3}$ is an annulus in $X(I_{1} \cup I_{3})$, but not in $X(J_{0,9,10})$ (it intersects $I_{2}$ and $I_{4}$).
If we twist along $A_{2,4}$, we destroy $A_{1,3}$.
However, it is not hard to see that we may remove a disk from $A_{1,3}$ so that the leftover surface misses a neighborhood of $A_{2,4}$, apply the twist in $A_{2,4}$, then glue a disk back to obtain another annulus we also call $A_{1,3}$.
We are using the fact that the regular neighborhood of $A_{2,4}$ is a solid torus and that the curve on the boundary which bounds a disk continues to bound a disk in $S^{3}$ after twisting.

Next, we observe that $I_{0} \cup I_{9}$ is the two component unlink.
Moreover, the above fillings do not change this.
That is $I_{0}(\vec{r}) \cup I_{9}(\vec{r})$ is still a two component unlink.
In fact, the basis for $H_{1}(\partial_{9} X(I_{0}(\vec{r}) \cup I_{9}(\vec{r})))$ is still $m_{9},l_{9}$ (that is, the new embedding into $S^{3}$ does not change the curve which bounds a meridian disk, nor the element which spans a Seifert surface).
The same is not true of the basis for $H_{1}(\partial_{0} X(I_{0}(\vec{r}) \cup I_{9}(\vec{r})))$, however $m_{0}$ still bounds a meridian disk.

We may therefore Dehn twist in the disk $D_{9}$ bounded by $I_{9}$ to obtain another embedding of $X(J_{0,10})$ into $S^{3}$.
Thus, the $(1,r_{5})$ filling on $\partial_{9}X(J_{0,10})$ for any $r_{5} \in {\mathbb Z}$, in addition to filling the other components according to (\ref{slopeseqn}), gives $M_{0,10}(\vec{r}) \cong S^{3}$ with $I_{0}(\vec{r})$ unknotted.
This proves parts 1 and 2 of the claim.

\begin{figure}[htb]
\begin{center}
\ \psfig{file=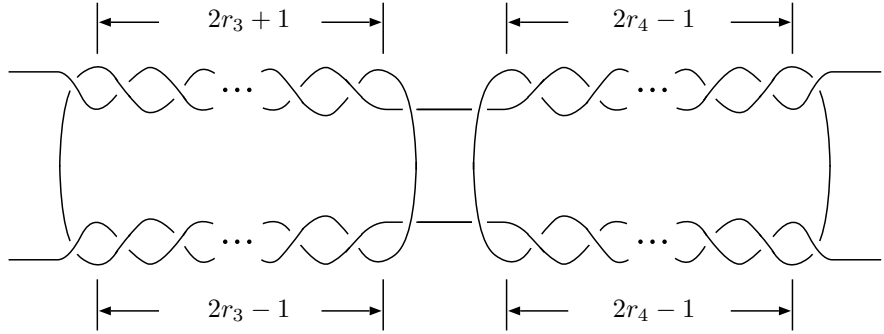,height=1.75truein}
\caption{Local picture of $U$ after filling (the numbers represent the number of crossings).}
\label{neari0}
\end{center}
  \setlength{\unitlength}{1in}
  \begin{picture}(0,0)(0,0)
    \put(1.78,2.22){$2r_{3}+1$}
    \put(1.78,.7){$2r_{3}-1$}
    \put(3.9,2.22){$2r_{4}-1$}
    \put(3.9,.7){$2r_{4}-1$}
  \end{picture}
\end{figure}

All that remains is to verify part 3 of the claim.
Consider a ball $U$ containing $I_{0}$, $I_{5}$, $I_{6}$, $I_{7}$, and $I_{8}$ as well as an arc of $I_{9}$ and $I_{10}$, but no other part of $J$.
The boundary of $U$ is indicated by a dotted circle in Figure \ref{gettintheknot}.
After filling the boundary components $\partial_{5}X(J_{0,...,4,9,10}) \cup \cdots \cup \partial_{8}X(J_{0,...,4,9,10})$, the neighborhood $U$ and the parts of $I_{0}(\vec{r})$, $I_{9}(\vec{r})$, and $I_{10}(\vec{r})$ in this neighborhood, are as in Figure \ref{neari0}.
The disk bounded by $I_{0}$ is contained in $U$, so $I_{10}$ intersects this disk only at points inside $U$.
Therefore, only the fillings on $\partial_{9}X(J_{0,10})$ may affect this intersection number (the annuli $A_{1,3}$ and $A_{2,4}$ are disjoint from $U$).
However, since $I_{9} \cup I_{10}$ is the unlink, the filling on $I_{9}$ also does not change this intersection number.
So, the linking number $\lk(I_{0}(\vec{r}),I_{10}(\vec{r}))$ (equivalently, this intersection number) is one.
This completes the proof of the claim. \hfill $\Box$\\

We are now finally in a position to describe the knots and link which prove Lemma \ref{knotlemma}.
The $p$-fold cyclic branched cover of $M_{0,10} \cong S^{3}$ (by part 1 of the claim) branched over $I_{0}(\vec{r})$ is $S^{3}$, since $I_{0}(\vec{r})$ is unknotted (by part 2 of the claim).
We let $K_{p}(\vec{r})$ denote the preimage of $I_{10}(\vec{r})$ under this covering, which is a knot by part 3 of the claim.

We can view the restriction of this branched cover to $X(K_{p}(\vec{r}))$ as a manifold covering of an orbifold
$$f: X(K_{p}(\vec{r})) \rightarrow X(I_{0}(\vec{r}) \cup I_{10}(\vec{r}))((p,0),\infty) = {\cal O}_{p}(\vec{r})$$
(the equality is (\ref{fillingeqn})).
We also view ${\cal O}_{p}$ as a suborbifold of ${\cal O}_{p}(\vec{r})$.
As such, its preimage under the covering is a submanifold $M_{p}(\vec{r}) \subset X(K_{p}(\vec{r}))$, for which $f$ restricts to a covering of ${\cal O}_{p}$.

We first note that, for fixed $p$, there are only finitely many homeomorphism types of $M_{p}(\vec{r})$.
This is because there are only finitely many $p$-fold covers of ${\cal O}_{p}$.
Second, we observe that $M_{p}(\vec{r}) = X(L_{p}(\vec{r}))$ for some link $L_{p}(\vec{r})$ in $S^{3}$, one component of which is $K_{p}(\vec{r})$; we take $L_{p}(\vec{r})$ to be the preimages of the cores of the filled in solid tori, union with $K_{p}(\vec{r})$, under the branched cover.

Fix any $p \geq 3$, and take any sequence $\{\vec{r}(j)\}_{j=1}^{\infty}$ for which $R < |r_{i}(j)| \rightarrow \infty$ as $j \rightarrow \infty$, with $\{ r_{i}(j) \}_{j=1}^{\infty}$ a sequence of distinct integers, for each $i = 1,...,5$, and so that $M_{p}(\vec{r}(j)) \cong M_{p}(\vec{r}(j'))$ for every $j,j'= 1,2,..$.

Set $L = L_{p}(\vec{r}(1))$ and $K_{j} = K_{p}(\vec{r}(j))$.
Then $X(L) = M_{p}(\vec{r}(1)) \cong M_{p}(\vec{r}(j))$, and hence each $X(K_{j})$ is obtained by Dehn filling all but one component of $X(L)$, for every $j$.
By construction, $X(L)$ and $X(K_{j})$ are all hyperbolic.
Moreover, $X(L)$ is a cover of ${\cal O}_{p}$, and thus contains the totally geodesic surface $F = f^{-1}(T_{p}^{*})$.
$F$ is a cyclic $p$-fold branched cover of a torus branched over two points, and thus has genus $p$.
This proves parts 1 and 3 of the lemma and most of part 2.
All that remains is to show that the slopes on each cusp of $X(L)$ are distinct.
To see this, we note that the filling slopes on the cusps of $X(L)$ are determined by lifting the filling slopes on the cusps of ${\cal O}_{p}$.
These latter are all distinct by hypothesis, and thus the former are as well. \hfill $\Box$

\subsection{Cusped totally geodesic surfaces} \label{knotcuspsection}

Here we describe a construction of hyperbolic knots in $S^{3}$ containing totally geodesic cusped surfaces in their complements.
In \cite{Adams3} and \cite{Adams4}, totally geodesic cusped surfaces such as these (and others) are constructed.
Moreover, in \cite{Adams4} the authors obtain some necessary conditions for a knot to contain such a surface as well as some uniqueness results for these types of surfaces.

\begin{example} \label{cuspedsurfaceexample} {\bf (compare \cite{Adams3} and \cite{Adams4})}
There exists hyperbolic knots in $S^{3}$ which contain embedded totally geodesic cusped surfaces in their complements.
\end{example}

\begin{figure}[htb]
\begin{center}
\ \psfig{file=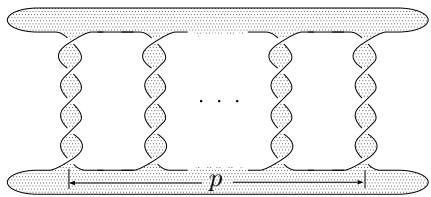,height=1truein}
\caption{A totally geodesic surface (for $p \geq 3$ and odd).}
\label{cuspedgeodesic}
\end{center}
  \setlength{\unitlength}{1in}
  \begin{picture}(0,0)(0,0)
    \put(2.97,.67){$p$}
  \end{picture}
\end{figure}

\begin{figure}[htb]
\begin{center}
\ \psfig{file=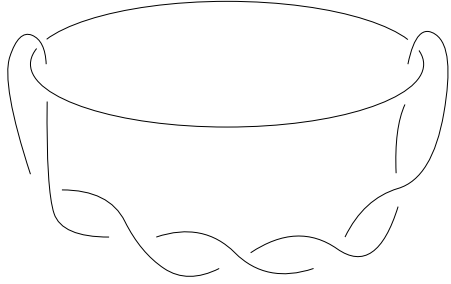,height=1.5truein}
\caption{A 2-component chain link.}
\label{duh}
\end{center}
\end{figure}

Begin with the 2-component chain link shown in Figure \ref{duh}, which was shown to be hyperbolic in \cite{NeumannReid}.
Let $K_{1}$ and $K_{2}$ denote the two components.
By Theorem \ref{hyperbolicdehnsurgery}, $X(K_{1} \cup K_{2})((p,0),\infty)$ is hyperbolic for $p$ sufficiently large, and as in \S \ref{linksection}, it suffices to take $p \geq 3$.
The $p$-fold cyclic branched cover of $S^{3}$ branched over $K_{1}$ is $S^{3}$ since $K_{1}$ is unknotted.
The complement of the preimage of $K_{2}$ under this branched covering is a link $K_{p} \subset S^{3}$ and we view the restriction of the branched cover to the complement of $K_{p}$ as a manifold covering of the orbifold
$$f_{p}:X(K_{p}) \rightarrow X(K_{1} \cup K_{2})((p,0),\infty)$$
If $p$ is odd, then $K_{p}$ is connected (since $\lk(K_{1},K_{2}) = 2$), which is to say, $K_{p}$ is a knot.
When $p \geq 3$, we thus obtain a hyperbolic knot $K_{p}$.

Since $K_{2}$ is also unknotted, it bounds a disk, and $K_{1}$ intersects this disk exactly twice.
Therefore, after $(p,0)$ surgery on $K_{1}$, this disk is a $(p,p,\infty)$ triangle orbifold, $T_{p}$, and hence is totally geodesic.
$f_{p}^{-1}(T_{p})$ is a totally geodesic cusped surface in $X(K_{p})$.

These knots are pictured in Figure \ref{cuspedgeodesic}.
The shaded Seifert surface is totally geodesic.

\section{Small curvature surfaces are acylindrical} \label{smallcurvesurface}

Let $F$ be a closed oriented surface in an oriented finite volume hyperbolic 3-manifold $M$.
As mentioned in the introduction, if $F$ is totally geodesic, then not only is it incompressible, but $M / F$ is acylindrical.
To see this, note that if $M/F$ contained an essential annulus, then $M/ F$ doubled along $F$ would contain an essential torus.
This is impossible since the doubled manifold is clearly hyperbolic.

An easy, unpublished result of Thurston states that to guarantee incompressibility, the requirement that $F$ is totally geodesic can be relaxed to a principal curvature bound.
More precisely, one has

\begin{theorem} [Thurston] \label{principal2incompressible}
If $F$ is a closed orientable surface in an orientable finite volume hyperbolic 3-manifold $M$ and has all principal curvatures less than $1$ in absolute value, then $F$ is incompressible.
Moreover, $F$ is quasi-Fuchsian.
\end{theorem}

For completeness, and since we will use it later, we give the proof.\\

\noindent
{\em Proof.}
By compactness, there is a global bound, ${\cal K} < 1$, for the absolute values of the principal curvatures, $\lambda_{1},\lambda_{2}$, of $F$.
It follows from Theorem \ref{gauss} that the (Gaussian) curvature $K(p) = K(T_{p}(F))$ of $F$ at $p$ satisfies
$$K(p) - (- 1) = \lambda_{1}(p)\lambda_{2}(p)$$
so that
$$|K(p) - (-1)| = | \lambda_{1}(p)\lambda_{2}(p) | \leq {\cal K}^{2} < 1$$
an hence $K(p)$ is negative for all $p$.

Let $\pi_{M} : {\mathbb H}^{3} \rightarrow M$ denote the universal covering of $M$, $\pi_{F} : \widetilde{F} \rightarrow F$ the universal covering of $F$, and $\widetilde{\pi}_{F} : \widetilde{F} \rightarrow {\mathbb H}^{3}$
some lift of $\pi_{F}$ composed with the inclusion of $F$ into $M$.

Because $F$, and hence $\widetilde{F}$, is negatively curved, any two points $p$ and $q$ in $\widetilde{F}$ are connected by a unique geodesic segment $\gamma$.
The geodesic curvature of $\widetilde{\pi}_{F} \circ \gamma$ for any $t$ satisfies $\kappa_{\widetilde{\pi}_{F} \circ \gamma}(t) \leq {\cal K}$ by (\ref{geodcurvboundeqn}) of \S \ref{Riemanniansection}, and $\widetilde{\pi}_{F} \circ \gamma$ is thus a $\frac{1}{\sqrt{1-{\cal K}^2}}$-quasi-geodesic, by Lemma \ref{quasicurvatureclose}.
In particular, this implies that $\widetilde{\pi}_{F}$ is a $\frac{1}{\sqrt{1-{\cal K}^{2}}}$-quasi-isometric embedding, and so $F$ is quasi-Fuchsian. \hfill $\Box$\\

Not only is the surface incompressible, but if the principal curvatures are small enough, we can (almost) recover acylindricity of the cut-open manifold.
In particular, we have

\begin{theorem} \label{curvature_acylindrical}
Given $g, r > 0$, there exists $\delta > 0$ with the following property.
If $F$ is an embedded, closed, orientable surface in an oriented finite volume hyperbolic 3-manifold $M$, with $genus(F) \leq g$, $injrad(F) \geq r$, and principal curvatures of $F$ bounded by $\delta$ in absolute value, then either $F$ bounds a twisted $I$-bundle in $M$, or $M/ F$ is acylindrical.
\end{theorem}

Although the proof is elementary, we divide it into a few lemmas for clarity.
Throughout, we write
$$\pi_{M}:{\mathbb H}^{3} \rightarrow M$$
to denote the universal covering of $M$.
We will write $B_{X}(p,R)$ to denote the closed ball in a metric space $X$ centered at a point $p$ with radius $R$.
$N_{X}(Y,R)$ will denote the closed $R$-neighborhood of a subset $Y$ of a metric space $X$.
When $X = {\mathbb H}^{3}$, we will simply write $B(p,R)$ and $N(Y,R)$ respectively.
If $(\widetilde{F},p) \subset {\mathbb H}^{3}$ is a pointed surface (i.e. a surface together with a point on the surface), we write
$${\mathbb H}^{2}(\widetilde{F},p)$$
to denote the unique (totally geodesic) hyperbolic plane in ${\mathbb H}^{3}$ tangent to $\widetilde{F}$ at the point $p$. 
We will write $d$ to denote the distance function on ${\mathbb H}^{3}$ and $d_{H}$ the corresponding Hausdorff distance.

We begin with

\begin{lemma} \label{close_stays_close}

Given $\epsilon > 0$ and $R > 0$, there exists $\delta_{0} > 0$ and $0 < \epsilon_{0} < \epsilon$, such that if
$$(\widetilde{F}_{1},p_{1}),(\widetilde{F}_{2},p_{2}) \subset {\mathbb H}^{3}$$
are disjoint, complete, embedded, pointed surfaces with principal curvatures bounded by $\delta_{0}$ in absolute value and $d(p_{1},p_{2}) < \epsilon_{0}$, then
$$d_{H}(B_{\widetilde{F}_{1}}(p_{1},R),B_{\widetilde{F}_{2}}(p_{2},R)) < \epsilon$$
\end{lemma}

\vspace{.5cm}

Roughly speaking, this lemma says that two disjoint surfaces in ${\mathbb H}^{3}$ with small principal curvature, which are close at some pair of points, must be close on large disks about those points.\\

\noindent
{\em Proof. }
Given $\eta > 0$ there exists $0< \mu < 1$, such that if
$$\gamma: [0,R+1] \rightarrow {\mathbb H}^{3}$$
is a unit speed path with $\kappa_{\gamma}(t) < \mu$, and if 
$$\sigma: [0,R+1] \rightarrow {\mathbb H}^{3}$$
is the unique geodesic with $\dot{\sigma}(0) = \dot{\gamma}(0)$, then for all $t \in [0,R+1]$,
$$d(\gamma(t),\sigma(t)) < \eta$$
This essentially follows from the proof of Lemma \ref{quasicurvatureclose}.

Suppose that $(\widetilde{F},p)$ is a complete, embedded, pointed surface in ${\mathbb H}^{3}$ with principal curvatures bounded by $\mu$ in absolute value.
Consider the exponential maps
$$\exp_{{\mathbb H}^{2}(\widetilde{F},p)},\exp_{(\widetilde{F},p)} : T_{p}(\widetilde{F}) \rightarrow {\mathbb H}^{3}$$
(this makes sense because ${\mathbb H}^{2}(\widetilde{F},p)$ and $\widetilde{F}$ are tangent at $p$).
By (\ref{geodcurvboundeqn}) of \S \ref{Riemanniansection}
$$d(\exp_{{\mathbb H}^{2}(\widetilde{F},p)}(v),\exp_{(\widetilde{F},p)}(v)) < \eta$$

In particular, note that if $\{ (\widetilde{F}_{n},p_{n}) \}_{n=1}^{\infty}$ is any sequence of pointed surfaces with principal curvatures approaching $0$ as $n \rightarrow \infty$ and if $\{ p_{n} \}_{n=1}^{\infty}$ has compact closure in ${\mathbb H}^{3}$, then after appropriately isometrically reparameterizing the domains (to ${\mathbb R}^{2}$ say) there is subsequence of $\{ \exp_{(\widetilde{F}_{n},p_{n})}|_{B_{{\mathbb R}^{2}}(0,R)} \}$ which converges uniformly by the Arzela-Ascoli Theorem (see e.g. \cite{Munkres}).
Moreover, the limit is easily seen to be (a reparameterization of) the exponential map restricted to a radius $R$ ball, with image a radius $R$ ball in a hyperbolic plane in ${\mathbb H}^{3}$.

Now suppose there is no $\epsilon_{0}$ and $\delta_{0}$ as in the statement of the theorem.
This implies that there exists a pair of sequences of pointed surfaces
$$\{ (\widetilde{F}_{1,n},p_{1,n}) \}_{n=1}^{\infty} \, \mbox{ and } \, \{ (\widetilde{F}_{2,n},p_{2,n}) \}_{n=1}^{\infty}$$
such that as $n \rightarrow \infty$, the principal curvatures of $\widetilde{F}_{i,n}$ approach $0$, for $i=1,2$, and $d(p_{1,n},p_{2,n}) \rightarrow 0$.
Moreover, $\widetilde{F}_{(1,n)} \cap \widetilde{F}_{2,n} = \emptyset$ and
$$d_{H}( B_{\widetilde{F}_{1,n}}(p_{1,n},R), B_{\widetilde{F}_{2,n}}(p_{2,n},R)) \geq \epsilon$$
for all $n \in {\mathbb Z}^{+}$.
By composing the embeddings with an isometry of ${\mathbb H}^{3}$, we can assume that $p_{1,n}$ is the same point for all $n$.
The above remarks imply that by passing to a subsequence, we may assume that $\{ \exp_{(\widetilde{F}_{i,n},p_{i,n})}|_{B_{{\mathbb R}^{2}}(0,R)} \}_{n=1}^{\infty}$ converges uniformly for each $i=1,2$ (after reparameterizing).

If the images of the two limit exponential maps are radius $R$ disks in distinct hyperbolic planes in ${\mathbb H}^{3}$, then they must non-trivially transversely intersect since the sequences of base points $\{p_{1,n}\}_{n=1}^{\infty}$ and $\{p_{2,n}\}_{n=1}^{\infty}$ must converge to a single point.
It then follows that $\widetilde{F}_{1,n}$ and $\widetilde{F}_{2,n}$ must intersect for sufficiently large $n$, which is a contradiction.
Therefore, the images of the two limit exponential maps are the same radius $R$ disk in a hyperbolic plane in ${\mathbb H}^{3}$.
This easily implies 
$$d_{H}( B_{\widetilde{F}_{1,n}}(p_{1,n},R), B_{\widetilde{F}_{2,n}}(p_{2,n},R)) < \epsilon$$
for sufficiently large $n$, which is also a contradiction.
It follows that the required $\epsilon_{0}$ and $\delta_{0}$ exists. \hfill $\Box$\\

\begin{lemma}  \label{3_close_separate}
There exists $\delta_{1} > 0$ and $\epsilon > 0$, such that if $\{ (\widetilde{F}_{i},p_{i}) \}_{i=1}^{3}$ are pairwise disjoint, complete, pointed surfaces in ${\mathbb H}^{3}$ with principal curvatures bounded by $\delta_{1}$ in absolute value and if $d(p_{1},p_{2}),d(p_{1},p_{3}) < \epsilon$, then one of the three surfaces separates ${\mathbb H}^{3}$ into two components, each of which contains exactly one of the remaining two surfaces.
\end{lemma}

\noindent
{\em Proof. }
We begin by noting that if $\widetilde{F}$ is a complete surface in ${\mathbb H}^{3}$ whose principal curvatures are bounded in absolute value by a constant less than $1$, then as in the proof of Theorem \ref{principal2incompressible}, the inclusion of $\widetilde{F}$ is a quasi-isometric embedding.
$\widetilde{F}$ is thus properly embedded, and so separates ${\mathbb H}^{3}$ into two components.

We now proceed in a fashion similar to that of the proof of Lemma \ref{close_stays_close}.
If no such $\delta_{1}$ or $\epsilon$ as in the statement of the lemma exist, then there must be three sequences of pointed surfaces
$$ \{ (\widetilde{F}_{i,n},p_{i,n}) \}_{n=1}^{\infty} \mbox{ for } i=1,2,3$$
with $\widetilde{F}_{i,n}$ having principal curvatures bounded in absolute value by $\frac{1}{n}$,and so that $d(p_{1,n},p_{2,n}),d(p_{1,n},p_{3,n}) < \frac{1}{n}$ and no one of the surfaces separates the other two.

As in the proof of Lemma \ref{close_stays_close}, by reparameterizing and passing to a subsequence if necessary, we may assume that the exponential maps for $B_{\widetilde{F}_{i,n}}(p_{i,n},10)$ converge uniformly to the exponential map for $B_{{\mathbb H}^{2}}(p,10)$ of some pointed hyperbolic plane $({\mathbb H}^{2},p)$ in ${\mathbb H}^{3}$ for each $i=1,2,3$ (the choice of radius $10$ here is arbitrary).
For sufficiently large $n$, the intersections $\widetilde{F}_{1,n} \cap B(p,9),\widetilde{F}_{2,n} \cap B(p,9),$ and $\widetilde{F}_{3,n} \cap B(p,9)$ are properly embedded disks in $B(p,9)$ (in fact, one can show that any complete surface $\widetilde{F}$ in ${\mathbb H}^{3}$ with principal curvatures bounded in absolute value by a constant less than $1$ must intersect a closed ball in a convex subset of $\widetilde{F}$).
As $n$ approaches infinity, the boundaries of these three disks converge uniformly to a great circle.
So for sufficiently large $n$, the three boundaries consists of three parallel loops on $\partial B(p,9)$.
One of these loops must separate $\partial B(p,9)$ into two components, each containing exactly one of the other two loops. 
The corresponding surface separates ${\mathbb H}^{3}$ into two components, each containing exactly one of the other two surfaces.
This contradiction proves the lemma. \hfill $\Box$\\

The next lemma says that an essential annulus in the manifold obtained by cutting open along a small principal curvature surface forces two components of the preimage of the surface in ${\mathbb H}^{3}$ to be close together.

\begin{lemma} \label{annulus_then_close}
Given $\epsilon_{0} > 0$, there exists a $\delta_{2} >0$, such that for any embedded, closed, oriented surface $F$ in an oriented finite volume hyperbolic 3-manifold $M$, if the principal curvatures of $F$ are bounded by $\delta_{2}$ in absolute value and $M/ F$ contains an essential annulus, then there are two pointed components $(\widetilde{F}_{0},p_{0}),(\widetilde{F}_{1},p_{1}) \subset \pi_{M}^{-1}(F)$ such that $d(p_{0},p_{1}) < \epsilon_{0}$.
\end{lemma}

\noindent
{\em Proof. }
By Lemma \ref{quasicurvatureclose}, there exists $\delta_{2} > 0$, such that if $\gamma : {\mathbb R} \rightarrow {\mathbb H}^{3}$ is a path with $\kappa_{\gamma}(t) < \delta_{2}$, then there is a unique geodesic $g_{\gamma} : {\mathbb R} \rightarrow {\mathbb H}^{3}$ which remains a bounded Hausdorff distance from $\gamma$, and moreover
$$d_{H} (\gamma,g_{\gamma}) < \frac{\epsilon_{0}}{2}$$

Now suppose that $F \subset M$ is as above, with principal curvatures bounded by $\delta_{2}$ in absolute value and
$$A: (S^{1} \times [0,1], \partial (S^{1} \times [0,1])) \rightarrow (M,F)$$
is an essential annulus with both boundary components on $F$.
We may homotope $A$ to guarantee that $A(\partial (S^{1} \times [0,1]))$ consists of two geodesics in $F$.

Next let
$$\widetilde{A}: {\mathbb R} \times [0,1] \rightarrow {\mathbb H}^{3}$$
be a lift of $A$.
Because $A$ is essential, $\widetilde{A}$ maps the boundary into two distinct components of $\pi^{-1}(F)$, which we call $\widetilde{F}_{0}$, and $\widetilde{F}_{1}$.
The maps of the boundary,
$$\gamma_{i} = \widetilde{A}|_{{\mathbb R} \times \{ i \}}: {\mathbb R} \rightarrow \widetilde{F}_{i}$$
for $i=0,1$, are geodesics in each $\widetilde{F}_{i}$ so that the geodesic curvatures satisfy $\kappa_{\gamma_{i}} < \delta_{2}$.
$\widetilde{A}$ defines a hyperbolic transformation (from the action of $\pi_{1}(M)$) having an axis $g$ and this transformation stabilizes each $\gamma_{i}$, $i=0,1$.
It follows that $g$ is the unique geodesic with $d_{H} (\gamma_{i},g) < \frac{\epsilon_{0}}{2}$ for $i=0,1$.
Therefore $d_{H} (\gamma_{0},\gamma_{1}) < \epsilon_{0}$.
In particular, this implies that there is a point on $\widetilde{F}_{0}$ closer than $\epsilon_{0}$ to a point on $\widetilde{F}_{1}$. \hfill $\Box$ \\

We now proceed with the proof of Theorem \ref{curvature_acylindrical}.\\

\noindent
{\em Proof. }
There exists a number $R > 0$, depending on $g$, $r$, and a fixed positive number less than $1$, say $\frac{1}{2}$, so that given any hyperbolic manifold $M$ and closed surface $F$ with $genus(F) \leq g$, $injrad(F) \geq r$, and  principal curvatures bounded in absolute value by $\frac{1}{2}$, then the diameter of $F$ is bounded by $\frac{R}{2}$.
This follows from the fact that the (Gaussian) curvature of $F$ is pinched between two negative constants (by Theorem \ref{gauss}) which gives an upper bound on the area of $F$ by the Gauss-Bonnet Theorem \cite{ONeill} and a lower bound on the area of an embedded disk of radius $r$ by an application of the Rauch Comparison \cite{DoCarmo} and Gauss-Bonnet Theorems, for example.

Let $\delta_{1},\epsilon$ be as in Lemma \ref{3_close_separate}, choose $\delta_{0},\epsilon_{0} > 0$ from Lemma \ref{close_stays_close}, based on $\epsilon$ and $R$, and choose $\delta_{2}$ from Lemma \ref{annulus_then_close}, based on $\epsilon_{0}$.
Now, set
$$\delta = \min \{ \frac{1}{2},\delta_{0},\delta_{1},\delta_{2} \}$$

Let $F$ be an embedded, closed, oriented surface in an oriented, finite volume, hyperbolic 3-manifold, with principal curvature bounded by $\delta$ in absolute value, $genus(F) \leq g$, and $injrad(F) \geq r$.
To prove the theorem, we suppose that $M/ F$ contains an essential annulus and prove that $F$ bounds a twisted $I$-bundle.

By Lemma \ref{annulus_then_close}, there are two components of $\pi_{M}^{-1} (F)$ separated by a distance less than $\epsilon_{0}$.
Choose two components $\widetilde{F}_{1},\widetilde{F}_{2} \subset \pi_{M}^{-1}(F)$ which are closest (by compactness of $F$, there exists a closest pair of components).
By Lemma \ref{close_stays_close}, there are points $p_{1} \in \widetilde{F}_{1}$ and $p_{2} \in \widetilde{F}_{2}$ such that
\begin{equation} \label{distance_bound}
d_{H}(B_{\widetilde{F}_{1}}(p_{1},R),B_{\widetilde{F}_{2}}(p_{2},R)) < \epsilon
\end{equation}

Since the diameter of $F$ is less than $\frac{R}{2}$, there are generators $\gamma_{1},...,\gamma_{n}$ for $stab_{\pi_{1}(M)}({\widetilde{F}_{1}}) \cong \pi_{1}(F)$, such that
$$\gamma_{j}(p_{1}) \in B_{\widetilde{F}_{1}}(p_{1},R)$$
for each $j = 1,...,n$.
By (\ref{distance_bound}), there are points $p_{2,1},...,p_{2,n} \in \widetilde{F}_{2}$ such that
$$d(\gamma_{j}(p_{1}),p_{2,j}) < \epsilon$$
for $j=1,...,n$.\\

\noindent
{\bf Claim. }
$\gamma_{j}(\widetilde{F}_{2}) = \widetilde{F}_{2}$ for each $j=1,...,n$.\\

\noindent
{\em Proof of claim. }
Suppose that there exists $1 \leq j \leq n$ such that $\gamma_{j}(\widetilde{F}_{2}) \neq \widetilde{F}_{2}$.
Then
$$(\widetilde{F}_{1},\gamma_{j}(p_{1})),(\widetilde{F}_{2},p_{2,j}),(\gamma_{j}(\widetilde{F}_{2}), \gamma_{j}(p_{2}))$$
are three disjoint pointed surfaces with $d(\gamma_{j}(p_{1}),p_{2,j}) < \epsilon$ and $d(\gamma_{j}(p_{1}),\gamma_{j}(p_{2})) = d(p_{1},p_{2}) < \epsilon$.
By Lemma \ref{3_close_separate}, one of the surfaces must separate ${\mathbb H}^{3}$ into two components, each of which contains exactly one of the other two surfaces.

Now, since $M$ and $F$ are orientable, $\gamma_{j}$ preserves the components of ${\mathbb H}^{3} \setminus \widetilde{F}_{1}$, and hence $\widetilde{F}_{2}$ and $\gamma_{j}(\widetilde{F}_{2})$ lie in the same component of ${\mathbb H}^{3} \setminus \widetilde{F}_{1}$.
If $\gamma_{j}(\widetilde{F}_{2})$ were the separating surface, then any path from $\widetilde{F}_{1}$ to $\widetilde{F}_{2}$ would have to pass through $\gamma_{j}(\widetilde{F}_{2})$.
This implies that $\gamma_{j}(\widetilde{F}_{2})$ is strictly closer to $\widetilde{F}_{1}$ than $\widetilde{F}_{2}$, which is impossible since $\widetilde{F}_{1}$ and $\widetilde{F}_{2}$ were chosen to be closest.
We now note that $\widetilde{F}_{1}$ and $\gamma_{j}(\widetilde{F}_{2})$ are also a pair of closest components (since $\gamma_{j}$ is an isometry), so by the same argument just given, $\widetilde{F}_{2}$ cannot separate these surfaces either.

Therefore no one of the surfaces can separate the other two, which is a contradiction.
This establishes the claim.\\

Let $\widetilde{X}$ denote the 3-manifold in ${\mathbb H}^{3}$ bounded by $\widetilde{F}_{1}$ and $\widetilde{F}_{2}$.
$\widetilde{X}$ is easily seen to be simply connected, irreducible, and invariant under $stab_{\pi_{1}(M)}(\widetilde{F}_{1})$ (since its boundary components are invariant by generators of $stab_{\pi_{1}(M)}(\widetilde{F}_{1})$).
The quotient,
$$X = \widetilde{X}/stab_{\pi_{1}(M)}(\widetilde{F}_{1})$$
is thus an irreducible 3-manifold homotopy equivalent to a surface.
It follows from standard 3-manifold topology (see e.g. \cite{Hempel} and \cite{Tnotes}) that $X \cong F \times [0,1]$.

Next, let $\gamma \in \pi_{1}(M)$ be such that
$$\gamma(\widetilde{F}_{1}) = \widetilde{F}_{2}$$
(which must exist since $\widetilde{F}_{1}$ and $\widetilde{F}_{2}$ are both components of $\pi_{M}^{-1}(F)$).
Either $\gamma(\widetilde{X}) = \widetilde{X}$ or $\gamma(\widetilde{X}) \cap \widetilde{X} = \widetilde{F}_{2}$.
However, if the latter held, then it must be the case that
$${\mathbb H}^{3} = \bigcup_{k \in {\mathbb Z}} \gamma^{k}(\widetilde{X})$$
and hence $M$ (or a two-fold cover) must fiber over the circle with fiber $F$.
By Theorem \ref{principal2incompressible}, $F$ is quasi-Fuchsian and hence cannot be a virtual fiber (see e.g. \cite{Tnotes}).

Therefore $\gamma(\widetilde{X}) = \widetilde{X}$, and again standard 3-manifold topology implies
$$X' = \widetilde{X} / < stab_{\pi_{1}(M)}(\widetilde{F}_{1}), \gamma >$$
is a twisted $I$-bundle and $X'$ embeds in $M$.
Therefore $F$ bounds a twisted $I$-bundle in $M$. \hfill $\Box$ \\

\section{Related questions} \label{relatedquestions}

The following two questions, related to Conjecture \ref{noknots}, seem to be interesting.

First, we note that links constructed in Section \ref{linksection}, as well as those constructed in \cite{MenascoReid}, have the property that the surface $F$ separates $S^{3}$ into two components (the two {\em sides} of the surface) with each side containing the same number of components of the link.
It is easy to construct examples for which there are a different number of components on each side of the surface (e.g. using techniques of \cite{Adams}).
However, all the constructions seem to require at least one component of the link on each side of the surface.
In particular, we have

\begin{question}
Do there exist hyperbolic links $L \subset S^{3}$ for which $X(L)$ contains an embedded, closed, connected, totally geodesic surface $F$ with $[F] = 0 $ in $H_{2}(X(L))$?
\end{question}

Next, we note that all known examples of one cusped hyperbolic manifolds with closed embedded totally geodesic surfaces contain non-peripheral homology.
In particular, these manifolds do not arise as knot complements in homology spheres.
We thus have

\begin{question} {\bf (Reid)}
Are there integral or rational homology spheres $M$ which contain hyperbolic knots with closed embedded totally geodesic surfaces in their complements?
\end{question}

\section{Appendix: Some computations in hyperbolic space} \label{appendix}

Here, as a convenience for the reader, we give a proof of\\

\noindent
{\bf Lemma \ref{quasicurvatureclose}}
{\em There exists a continuous non-negative function $f$ on the interval $[0,1)$, so that $f(0)=0$, having the following property.
Suppose $\gamma : [a,b] \rightarrow {\mathbb H}^{n}$ is a unit speed path whose geodesic curvature satisfies
$$\kappa_{\gamma}(t) \leq {\cal K}$$
for all $t \in [a,b]$, where $0 \leq {\cal K} < 1$ is some constant.
Then $\gamma$ is a $\frac{1}{\sqrt{1-{\cal K}^{2}}}$-quasi-geodesic.
Moreover, if $g_{\gamma}:[a,b] \rightarrow {\mathbb H}^{n}$ is the unique geodesic connecting the endpoints $\gamma(a)$ and $\gamma(b)$, then the Hausdorff distance between the image of $\gamma$ and $g_{\gamma}$ is no more than $f({\cal K})$.}\\

The first part of the proof we give here is a computational version of the proof Corollary 8.9.3 of \cite{Tnotes}.
The second part follows from the proof of Proposition 5.9.2 along with an application of the Arzela-Ascoli Theorem.\\

\noindent
{\em Proof.}
We use the hyperboloid model of hyperbolic space $H^{n} \subset {\mathbb E}^{n,1}$ (see \cite{Ratcliffe} for a complete description of this model).
We let $\langle,\rangle$ denote the form on ${\mathbb E}^{n,1}$, as well as its restriction to $T(H^{n})$.
Let $X,Y \in {\frak X} (H^{n})$, with $Y = Y^{i} \partial_{i}$, where $\partial_{i} = \frac{\partial}{\partial x^{i}}$.

The Levi-Civita connection on $H^{n}$ is then given by
$$\nabla_{X} Y (p)
= X(Y^{i})(p) \partial_{i} + \langle X(Y^{i})(p) \partial_{i}, p\rangle p$$
Here we view $p$ simultaneously as a point in $H^{n}$ and a vector in ${\mathbb E}^{n,1}$.
It follows that if $\gamma : I \rightarrow H^{n}$ is a path, and $V = V^{i} \partial_{i}$ is a vector field along $\gamma$, that the covariant derivative of $V$ is given by
$$\frac{DV}{dt}(t)
= \dot{V}(t) + \langle\dot{V}(t),\gamma(t)\rangle \gamma(t)$$
In particular, suppose $\gamma$ is a unit speed path, then we have
$$0 = \frac{d}{dt}(1) = \frac{d}{dt}\langle\dot{\gamma}(t),\dot{\gamma}(t)\rangle = \langle\ddot{\gamma}(t),\dot{\gamma}(t)\rangle + \langle\dot{\gamma}(t),\ddot{\gamma}(t)\rangle$$
and
$$0 = \frac{d}{dt}\langle\dot{\gamma}(t),\dot{\gamma}(t)\rangle = \langle\frac{D \dot{\gamma}}{dt}(t),\dot{\gamma}(t)\rangle + \langle\dot{\gamma}(t), \frac{D \dot{\gamma}(t)}{dt}(t)\rangle$$
So that
$$\langle\ddot{\gamma}(t),\dot{\gamma}(t)\rangle = 0 = \langle \frac{D \dot{\gamma}}{dt}(t),\dot{\gamma}(t)\rangle$$
We also see that
$$0 = \frac{d}{dt} (0) = \frac{d}{dt} \langle\dot{\gamma}(t),\gamma(t)\rangle = \langle\ddot{\gamma}(t),\gamma(t)\rangle + \langle\dot{\gamma}(t),\dot{\gamma}(t)\rangle$$
So that
$$\langle\ddot{\gamma}(t),\gamma(t)\rangle = -1$$
The geodesic curvature of $\gamma$ is given by $\kappa(t) = \| \frac{D \dot{\gamma}}{dt}(t) \|$, so that
\begin{align}
\|\ddot{\gamma}(t)\|^{2} & =
\| \ddot{\gamma}(t) + \langle\ddot{\gamma}(t),\gamma(t)\rangle \gamma(t) - \langle\ddot{\gamma}(t),\gamma(t)\rangle \gamma(t) \|^{2}
= \| \frac{D \dot{\gamma}}{dt}(t) + \gamma(t) \|^{2} \notag \\
 & = \| \frac{D \dot{\gamma}}{dt}(t) \|^{2} + \|\gamma(t)\|^{2}
= \kappa(t)^{2} - 1 \notag
\end{align}
Where the second to last equality holds because $\frac{D \dot{\gamma}}{dt}(t) \bot \gamma(t)$.

We now claim that if $\gamma$ is a unit speed path with $\kappa(t) < 1$, then the family of hyperbolic hyperplanes $P_{t}$ through $\gamma(t)$ and orthogonal to $\dot{\gamma}(t)$ are all disjoint (recall that $P_{t}$ is the intersection with $H^{n}$ of the linear subspace $(\dot{\gamma}(t))^{\perp}$).
To see this, note that for $t_{0}$ and $t_{0} + t$ in the domain of definition of $\gamma$, and $t > 0$, $P_{t_{0}}$ and $P_{t_{0}+t}$ are disjoint if the dual vectors $\dot{\gamma}(t_{0})$ and $\dot{\gamma}(t_{0}+t)$ span a subspace of signature $(1,1)$.
For then, $\langle,\rangle$ restricted to $(span\{\dot{\gamma}(t_{0}),\dot{\gamma}(t_{0}+t) \})^{\perp} = (\dot{\gamma}(t_{0}))^{\perp} \cap (\dot{\gamma}(t_{0}+t))^{\perp}$ is positive definite and so is disjoint from $H^{n}$.
The subspace has signature $(1,1)$ if there is a vector in their span with negative norm squared.
Now, since
$$\ddot{\gamma}(t_{0}) = \lim_{t \rightarrow 0} \frac{\dot{\gamma}(t_{0}+t) - \dot{\gamma}(t_{0})}{t}$$
we see that for $t$ sufficiently small, we have that $\dot{\gamma}(t_{0}+t) - \dot{\gamma}(t_{0})$ has negative norm squared since, by hypothesis, $\ddot{\gamma}(t_{0})$ does (because $\| \ddot{\gamma}(t) \|^{2} = \kappa(t)^{2} - 1$).

Since the hyperbolic hyperplanes $P_{t}$ are progressing in the direction of $\gamma$, and since locally they are being moved off themselves, we see that all the $P_{t}$'s are disjoint, proving the claim.

Next, we claim that at $t_{0}$, the hyperplanes $P_{t_{0}+t}$ are progressing at a rate of $\sqrt{1 - \kappa(t_{0})^{2}}$.
Specifically, we mean that
$$\lim_{t \rightarrow 0^{+}} \frac{d(P_{t_{0}},P_{t_{0}+t})}{ t } = \sqrt{1- \kappa(t_{0})^{2}}$$

To see this, we set $\delta(t) = d(P_{t_{0}},P_{t_{0}+t})$ and note that
$$ cosh(\delta(t)) = \langle\dot{\gamma}(t_{0}),\dot{\gamma}(t_{0}+t)\rangle$$
and that
\begin{align}
1 - \kappa(t_{0})^{2} & =  - \| \ddot{\gamma}(t_{0}) \|^{2}  =  - \lim_{t \rightarrow 0^{+}} \| \frac{\dot{\gamma}(t_{0}+t) - \dot{\gamma}(t_{0})}{t} \|^{2} \notag \\
 & =  - \lim_{t \rightarrow 0^{+}} \frac{\| \dot{\gamma}(t_{0}+t) \|^{2} - 2 \langle\dot{\gamma}(t_{0}+t),\dot{\gamma}(t_{0})\rangle + \| \dot{\gamma}(t_{0}) \|^{2}}{t^{2}} \label{tryitouteqn}\\
 & =  \lim_{t \rightarrow 0^{+}} \frac{2 \langle \dot{\gamma}(t_{0}+t),\dot{\gamma}(t_{0}) \rangle - 2}{t^{2}}
 = \lim_{t \rightarrow 0^{+}} \frac{ 2 (cosh(\delta(t)) - 1)  }{t^{2}} \notag \end{align}
We also see that since $\lim_{t \rightarrow 0^{+}} \delta(t) = 0$ (and since $\delta(t) > 0$ for $t > 0$), we can apply L'H$\hat{o}$spital's rule twice to obtain
$$ \lim_{t \rightarrow 0^{+}} \frac{\delta(t)^{2}}{2 ( cosh(\delta(t)) - 1)}
= \lim_{t \rightarrow 0^{+}} \frac{2\delta(t)}{2 sinh(\delta(t))}
= \lim_{t \rightarrow 0^{+}} \frac{1}{cosh(\delta(t))} = 1$$
Combining this with (\ref{tryitouteqn}), we obtain
\begin{align} 1 - \kappa(t_{0})^{2}
& = \lim_{t \rightarrow 0^{+}} \frac{ 2 ( cosh (\delta (t)) - 1 ) }{t^{2}}
= \lim_{t \rightarrow 0^{+}} \frac{ 2 ( cosh (\delta (t)) - 1 ) }{t^{2}} \frac{\delta(t)^{2}}{2 ( cosh(\delta(t)) - 1 ) } \notag \\
 & = \lim_{t \rightarrow 0^{+}} \frac{\delta(t)^{2}}{t^{2}}
= \lim_{t \rightarrow 0^{+}} \left( \frac{\delta(t)}{t} \right)^{2} \notag \end{align}
Taking square roots of the first and last term proves the claim.

We can therefore define the total displacement of the hyperplane $P_{t}$ from $P_{t_{0}}$ where $t > t_{0}$ as
$$\Delta(t) = \Delta_{t_{0}}(t) = \int_{t_{0}}^{t} \sqrt{1 - \kappa(s)^{2}} ds$$
If $t_{0} < t_{1} < ... < t_{n} = t$ is any partition of $[t_{0},t]$, then the quantity
$$\sum_{i=1}^{n} d(P_{t_{i-1}},P_{t_{i}})$$
approximates $\Delta(t)$ {\em from above}; any refinement of the partition does not increase the sum.
In particular, $\Delta(t) \leq d(P_{t_{0}},P_{t}) \leq d(\gamma(t_{0}),\gamma(t))$.
It follows that if $0 < \frac{1}{\lambda} \leq \sqrt{ 1 - \kappa(t)^{2} }$, for all $t$, then $\gamma$ is a $\lambda$-quasi-geodesic.

That is, if
$$\kappa(t) \leq {\cal K}$$
for all $t \in [a,b]$, then $\gamma$ is a $\frac{1}{\sqrt{1-{\cal K}^{2}}}$-quasi-geodesic.
This proves the first part of the lemma.\\

In the proof of Proposition 5.9.2 of \cite{Tnotes}, Thurston shows that a $\lambda$-quasi-geodesic segment remains a bounded distance from the geodesic with the same endpoints.
Moreover, he shows that this bound on the distance depends only on $\lambda$.
It is clear that if $\lambda_{0} < \lambda_{1}$ then the best bound does not increase.
However, the proof given does not imply that as $\lambda$ approaches $1$, then the best bound approaches $0$.
This is precisely what is needed to complete the proof.

We suppose therefore, that this best bound does not approach $0$, and arrive at a contradiction.
So, there exists $\epsilon>0$ and a sequence $\{ \gamma_{j} : [a_{j},b_{j}] \rightarrow {\mathbb H}^{n} \}_{j=1}^{\infty}$ so that each $\gamma_{j}$ is a unit speed $\lambda_{j}$-quasi-geodesic, $d_{H}(\gamma_{j}, g_{j}) > \epsilon$ where $g_{j}$ is the unique geodesic with the same endpoints as $\gamma_{j}$, and $\{ \lambda_{j} \}_{j=1}^{\infty}$ decreases to 1.
Note that because each $\lambda_{j}$ is bounded above by $\lambda_{1}$, it must be that $d_{H}(\gamma_{j},g_{j}) < R$ for some $R > 0$.

After reparameterizing each $\gamma_{j}$ (keeping it unit speed, but ensuring that $0 \in (a_{j},b_{j})$) and composing with an isometry of ${\mathbb H}^{n}$, we may assume that there is a geodesic line $g$ in ${\mathbb H}^{n}$ containing a point $x_{0}$, such that $g_{j} \subset g$ and
$$\epsilon < d(\gamma_{j}(0),g) = d(\gamma_{j}(0),x_{0})< R$$

By passing to a subsequence, we may further assume that as $j \rightarrow \infty$, $a_{j} \rightarrow a$ and $b_{j} \rightarrow b$ (with one or both of $a$ or $b$ possibly being infinite).

We now extend each $\gamma_{j}$ to a map on all of ${\mathbb R}$, by $\gamma_{j}(t) = \gamma_{j}(b_{j})$ for $t \in [b_{j},\infty)$ and $\gamma_{j}(t) = \gamma_{j}(t) = \gamma_{j}(a_{j})$ for $t \in (-\infty,a_{j}]$.
The set $\{ \gamma_{j} \}$ clearly forms a normal family, and so by passing to a subsequence, we may assume that the $\{ \gamma_{j} \}_{j=1}^{\infty}$ converges uniformly on compact subsets, by the Arzela-Ascoli Theorem.

Now we note that the limit $\gamma$, restricted to $(a,b)$ must be a $1$-quasi-geodesic, and so a geodesic.
Moreover, the ends must limit on $g$ (or the boundary of $g$ at infinity).
This implies that $\gamma \subset g$.
This is a contradiction since $\gamma_{j}(0)$ remains a distance at least $\epsilon$ from $g$, and hence so does $\gamma(0)$. \hfill $\Box$

\noindent
Address:\\
Department of Mathematics\\
Columbia University\\
2990 Broadway MC 4448\\
New York, NY 10027-6902\\
Phone: (512) 854-2431\\
email: clein@math.columbia.edu\\

\end{document}